\documentclass{amsart}
\usepackage[usenames,dvipsnames,svgnames,table]{xcolor}
\usepackage[margin=1.25in]{geometry}
\usepackage{amsfonts, amssymb}
\usepackage{enumerate}
\usepackage{graphicx}
\usepackage{amsmath}
\usepackage{enumitem}
\usepackage[normalem]{ulem} 
\usepackage{tikz}

\newtheorem{theorem}{Theorem}[section]

\newtheorem{claim}{Claim}

\newtheorem{cor}[theorem]{Corollary}

\newtheorem{prop}[theorem]{Proposition}

\newtheorem{lemma}[theorem]{Lemma}

\theoremstyle{definition}
\newtheorem{defn}[theorem]{Definition}

\newtheorem{remark}[theorem]{Remark}

\begin{document}
	\title{ Block Generalized Inverses in Skew Corner Pairs }
	\author{Feroz Siddique}
	\address{Department of Mathematics, University of Wisconsin Eau Claire, WI-54701, USA}
	\email{siddiquf@uwec.edu}
	
	\begin{abstract}
We view generalized inverses from another perspective, through invertibility in the skew corner bimodules $eRf$ attached to a pair of distinct idempotents $e, f$ of a ring $R$. Many generalized inverses are defined by a system of equations together with conditions on ranges and kernels, which makes them harder to visualize as inverses in the traditional sense. In the language of invertible blocks, however, they become genuine  one or two-sided inverses of a suitable compression of the element. We will  systematically show how block invertibility in a skew corner $eRf$ can be used to  visualize Drazin's $(b,c)$-inverse, which contains the Moore--Penrose, Bott-Duffin $(e,f)$ inverses amongst others. On this foundation we develop a way to produce units of the corner ring from complementary blocks, and transport stable range one across two-sided Peirce decompositions for distinct idempotents.  
	\end{abstract}
	
			\maketitle
		\noindent  All our rings are associative rings with identity and all our modules are left modules, unless specifically mentioned. $U(R)$ denotes the set of all units of $R$,  $J(R)$ denotes the Jacobson radical, $\operatorname{T}(a)$ denotes the set of all outer or pseudo-inverses of $a$, and $\operatorname{Idem}(R)$ denotes the set of  idempotents of the ring $R$. By $regular$ we mean $von$ $Neumann$ regular. By $\operatorname{sr}(R)$ we denote the stable range of the ring $R;$ and by $\operatorname{sr}_R(a)$, we denote the (left) stable range of the element $a$ in the ring $R$. \\


		An element $a$ in the skew corner $eRf$ is invertible if there is  some $a'\in fRe$ which satisfies $aa'=e$ and $a'a=f$. Such elements appear in the literature in several guises.  The idempotents $e$ and $f$ are algebraically \emph{Murray–von Neumann } equivalent precisely when $eRf$ contains an invertible element. This equivalence plays a fundamental role in nonstable algebraic $K$-theory. Two idempotents in $ M_n(R)$ are equivalent exactly when the corresponding modules are isomorphic. Then the equivalence classes of idempotents in $M_\infty(R)=\bigcup_{n\ge 1}M_n(R)$, with addition induced by block-diagonal sums, form the monoid $V(R)$ of isomorphism classes of finitely generated projective right $R$-modules \cite{kap2,murray}. From a semigroup-theoretic perspective, Green's relations~\cite{green} compare elements of $R$ through the one-sided ideals they generate: $x\,\mathcal{R}\,y$ if $xR=yR$, $x\,\mathcal{L}\,y$ if $Rx=Ry$, and $\mathcal{H}=\mathcal{R}\cap\mathcal{L}$. An element $a$ is an invertible block of $eRf$ exactly when $aR=eR$ and $Ra=Rf$. Here the invertible blocks of $eRf$ form a single $\mathcal{H}$-class (when there are any), and $a\mapsto a'$ maps them bijectively onto the invertible blocks of $fRe$.
		In the sense of Loos \cite{loos}, the pair $(eRf,fRe)$ is called an \emph{associative pair}, modelled on pairs of rectangular matrices, and invertible blocks are invertible elements of this pair. 
	Skew corners also enter the theory of stable range through Ara and Goodearl's \cite{goodearl1}  reduction of problems of higher stable rank to stable rank one for skew corners, and the $QB$-rings \cite{ara} which rest on quasi-invertible elements of skew corners, that is, elements $x\in pRq$ with $(p-xy)\perp(q-yx)$ for some $y\in qRp$ \cite{ara}. \\
	In the following sections we will  systematically show how block invertibility in a skew corner $eRf$ can be used to  visualize Drazin's $(b,c)$-inverse, which contains the Moore--Penrose, Drazin, group, core and dual core inverses, the Bott--Duffin inverse and the inverse along an element as special cases. We will  envison the stable range condition using our notion of one-sided block invertibility. As applications we  will provide key  results in the case of Bott–Duffin inverses, von Neumann regularity, unit-regularity, stable range and special cleanness.  For special-clean element, using block invertibility we provide a characterization of when every unit-regular element  is special clean, strengthening the theorem of Camillo and Khurana locally to elements.
		
		
		 

		\section{Introduction}
	\medskip 
	Following Ara and Goodearl \cite{goodearl1}, for idempotents $e,f\in R$  the subset $eRf$, naturally regarded as an $(eRe,fRf)$-bimodule, is called the \emph{skew corner} (also \emph{rectangular corner}) of $R$ determined by the pair $(e,f)$. Likewise $fRe$, an $(fRf,eRe)$-bimodule, is the \emph{skew corner} determined by $(f,e)$. When $e=f$, the skew corner $eRe$ is the classical corner ring of $R$, a ring with identity $e$. Every $a\in R$, compressed on either side by $e$ and $f$, produces the element $eaf\in eRf$, which we call the \emph{compression} or \emph{rectangular block} (or simply the \emph{block}) of $a$ in $eRf$.\\
	Before analyzing this relation, we fix some terminology for the subsets in $R$ in which our elements are associated with. 
		\begin{defn}
		Let	$e,f\in R$ be arbitrary idempotents. The skew-corner pairs $eRf$ and $fRe$ are called \emph{associate skew-corners}. A compression block $x\in eRf$ is called \emph{block invertible} in $eRf$ if there exists a compression block $y$ in the associate skew-corner $fRe$ such that
			\[
			xy=e \qquad\text{and}\qquad yx=f .
			\]
			In this case, $x$ and $y$ are called \emph{block inverses} of each other.
			\end{defn}
\noindent Kaplansky \cite{kap1} called idempotents $e$ and $f$ \emph{equivalent}, written $e\sim f$, if there exists a pair of block inverses $x\in eRf$ and $y\in fRe$. Equivalent idempotents are also called \emph{Murray--von Neumann equivalent} or \emph{Kaplansky equivalent}. Block invertibility is positional: whether an element is block invertible depends on the skew-corner in which it is considered and its inverse in the corresponding associate corner. For this reason, we sometimes record the skew-corners explicitly and speak of a block invertible pair $(x,y)$ in $(eRf,fRe)$.
This notion of invertible blocks in the skew corners  also appears as \emph{$(e,f)$-inverse}, or \emph{$(e,f)$-invertibility}; see Raki\'c and Djordjevi\'c \cite{rakic} and Mary \cite{mary1}.
	
	We  use the notation $e \overset{x}{\sim} f$ to  specify the invertible block $x= exf \in eRf.$ In this language, Kaplansky's definition reads: $e  \overset{x}{\sim} f$ if and only if the skew corner $eRf$ contains an invertible block $x$. Invertible blocks have exactly one unique block inverse. 

\noindent Representing an idempotent equivalence $e \sim f$ with $e \overset{x}{\sim} f$ allows us to use the chosen invertible block $x$ to serve as a  basis for the entire skew-corner in the form $eRf = x(fRf)$. This is because $e \overset{x}{\sim} f \implies \exists y \in fRe$ such that $xy = e $ and $yx = f.$ Then for any element $s \in eRf$ we have $ys \in fRf$ and \[ s = es = (xy)s= x(ys)\] 
This representation is unique because if  $xt_1=xt_2$ with $t_1,t_2\in fRf$, then \[ t_1=ft_1=yxt_1=yxt_2=ft_2=t_2 \] 
This gives us our bijection $t\mapsto xt$ from $fRf\to eRf$ with inverse $s\mapsto ys.$  
Including the invertible block in the equivalence notion also allows us to obtain the invertible block in a composition  \emph{i.e.}
\[e \overset{x}{\sim} f \overset{u}{\sim} g \implies e  \overset{xu}{\sim} g \]
Another observation  
is that we can represent the set of all invertible blocks in $eRf$ as either $x\operatorname{U}(fRf)$ or $\operatorname{U}(eRe)x$. 
We therefore have a one-to-one correspondence from  $U(fRf)$ or $U(eRe)$ onto the invertible blocks of $eRf$. 

 \noindent Continuing in this way, it is natural to ask what else survives of the above notions when the skew corner $eRf$ contains a one-sided inverse.




\begin{defn}
	Let $e,f\in R$ be idempotents and let $x\in eRf$. An element $y$ of the associate skew-corner $fRe$ is called a \emph{right block inverse} of $x$ if $xy=e$, and a \emph{left block inverse} of $x$ if $yx=f$. We say that $x$ is \emph{right} (resp.\ \emph{left}) \emph{block invertible} in $eRf$ if it admits a right (resp.\ left) block inverse in the associate skew corner $fRe$.
\end{defn}
Note that $x$ is block invertible  precisely when it is both left and right block invertible. In this case if $e \overset{x}{\sim} f$, then every right (resp. left) invertible block of $eRf$ is exactly $x$ times the right (resp. left)  invertible block elements of $fRf$. The same works for $eRe$ but by post-multiplication by $x$. In particular,we can see that when $e\sim f$, every right  or left  invertible block of $eRf$ $(fRe)$ is invertible if and only if either $fRf$ or $eRe$ is directly finite. This has interesting consequences that allows us to answer  when left $(b,c)$ invertibility is equivalent to right $(b,c)$ invertibility \cite[pg. ~71]{drazin2}, \cite{wang}. We shall see that  later in the next section when we visualize Drazin's $(b,c)$ inverse as a traditonal left-right invertibility notion. \\

Now the relation symbol captured by one-sided block invertibility is \emph{subequivalence}, for which we use the  notation below.
First recall that the  \emph{minus partial order} was defined by Hartwig\cite{hartwig1} and Nambooripad \cite{namboor}  on von Neumann regular elements $a, b \in R$ and given below :
\[a \leq^{-} b \text{ if and only if } a = eb = bf \text{ for  idempotents } e^2=e, f^2 = f \in R  \]


\noindent Then  the set  of idempotents of $R,$ $\operatorname{Idem}(R)$ is partially ordered by the minus partial order as follows: 
\[
p\leq^{-} q \iff p = pq = qp .
\]
\begin{defn}[Ara--Goodearl \cite{goodearl1}]\label{def:subequiv}
	Let $e,f\in R$ be idempotents. We say $e$ is \emph{subequivalent} to $f$,
	written $e\lesssim f$, if there exists an idempotent $g\in R$ such that
	\[
	e\sim g \qquad\text{and}\qquad g\leq^{-} f .
	\]
\end{defn}

\noindent In this case  $eR$ is isomorphic to a direct summand of $fR$ and this occurs precisely when  there exists elements  $x \in eRf$ and $y \in fRe$ such  that $xy = e.$ \\ 

\noindent We will  use $e \; {}^{x}{\lesssim} \; f$ to denote that the compression block $exf \in eRf$ being right invertible (\emph{i.e.} there exists a $w \in fRe$ such that $(exf)w = e$). Dually, for $y\in R$ we write $e\lesssim^{y}f$ if the compression $fye\in fRe$ is left invertible, that is, $z(fye)=e$ for some $z\in eRf$. The two notations describe the same relation $e\lesssim f$ from its two sides \emph{i.e.}  if $x\in eRf$ and $y\in fRe$ satisfy $xy=e$, then
\[
e\;{}^{x}{\lesssim}\;f\qquad\text{and}\qquad e\lesssim^{y}f
\]
both hold, and we record the relation with whichever of the two blocks is at hand. In each notation the block is written next to its left identity: $x=ex$ beside $e$, and $y=fy$ beside $f$.\\
We see that the relation $\lesssim$ is a preorder on $\operatorname{Idem}(R)$: it is
reflexive $(p \lesssim  p)$ and transitive. Indeed if $p \lesssim q$ and $ q \lesssim s$ for idempotents $p,q,s \in \operatorname{Idem}(R)$ then there exists elements $x_1 \in pRq, y_1 \in qRp$ and $x_2 \in qRs$ and $y_2 \in sRq$ such that $x_1y_1 = p$ and $x_2y_2 = q$. Consequently, 
\[(x_1x_2)(y_2y_1) = x_1(x_2y_2)y_1 = x_1(q)y_1= x_1y_1= p\]
where $(x_1x_2)\in pRs$ and $(y_2y_1) \in sRp$ and therefore, $ p \lesssim s$.\\
\medskip 
It is known \cite[Ex. ~ 2.2]{rizvi} that the relation $e\lesssim f\lesssim e$ need not give $e\sim f.$ However, if  the left and right block inverses of $x$ coincide then it is true  \emph{i.e.},
\[  e\;{}^{x}{\lesssim}\;f\ \text{ and }\ f\lesssim^{x}e\iff e\overset{x}{\sim}f \]
In general, a left or right  invertible  block in a skew corner need not be block invertible, even under stable range one assumption of the ambient ring. Take the following example of $R=M_2(k)$ with $f=I_2$ and $e=\begin{pmatrix}1&0\\0&0\end{pmatrix} $ which has stable range one.
Then $ e  \lesssim^{{\scriptscriptstyle fe}}  f $ since $(ef)(fe) = e$ but  $f \; {}^{\scriptscriptstyle r}{\lesssim} \; e$ is impossible for any $r \in fRe$ as that would force $rs = I_2$ for some $s \in eRf$. We will see in Sec. 3 that it requires a stronger form of each relation to force subequivalence of the idempotents into equivalence. \\ 

\noindent The following lemma shows that replacing either idempotent in $e\lesssim f$ by an equivalent one preserves the subequivalence. In particular, $\lesssim$ induces a well-defined preorder on the equivalence classes $[e]$ of idempotents, given by $[e]\le[f]$ if $e\lesssim f$.
\begin{lemma}{\label{chotabhim}}
	For idempotents $e,f,g \in \operatorname{Idem}(R)$,  the following hold. 
	\begin{enumerate}
		\item  $e \lesssim f  \sim g \implies e \lesssim g$ 
		\item $g \sim e \lesssim f \implies  g \lesssim f$
	\end{enumerate}
\end{lemma}
\begin{proof}
	$(1):$ Let $e,f \in \operatorname{Idem}(R)$ such that 
	\[ e \lesssim f\]
	Then by definition, there exists an idempotent $ p \in \operatorname{Idem}(R)$ such that $p \sim e$ and we have $ p \leq f \iff  p = fp = pf$. \\
	Now $p \sim e \implies $ there exists blocks $a\in pRe$ and $b\in eRp$ with \[ab=p \,\,  \text{ and }  \,\, 
	ba=e\]
	Now $ f\sim g \implies $there exists blocks $x\in fRg$ and $y\in gRf$ with \[xy=f \,\,  \text{ and }  \,\, 
	yx=g\] 
	Now set $ u=bx\in eRg, \qquad v=ya\in gRe. $ Then, 
	\[uv\;=\; b(xy)a = bp(xy)pa = b(pfp)a =  ba = e \]
	Therefore  $e\lesssim g$. \\
	\noindent $(2):$ The proof works similarly as above. \\ 
\end{proof}




\section{ Block Invertibility and Generalized Inverses}

\medskip 
Many generalized inverses are defined by a system of equations together with conditions on ranges and kernels. In the language of invertible blocks, where an invertible block is the ordinary notion of inverse read inside a skew corner, they become genuine two-sided inverses of a suitable compression of the element. We will first  illustrate this case with Drazin's $(b,c)$-inverse~\cite{drazin1} for arbitrary $b$ and $c.$  \\
\noindent Recall the following definition of $(b,c)$ inverse of an element.
\begin{defn}[Drazin \cite{drazin1,drazin2}]\label{bcinverse}
	Let $a,b,c\in R$. An element $y\in R$ is a \emph{$(b,c)$-inverse} of $a$ if
	\[
	y\in bRy\cap yRc,\qquad yab=b,\qquad cay=c.
	\]
	Such a $y$ is unique when it exists,  and  is usually denoted as $a^{(b,c)}$. \\
	$y$ is a \emph{left $(b,c)$-inverse} of $a$ if $Ry\subseteq Rc$ and $yab=b$, and a \emph{right $(b,c)$-inverse} of $a$ if $yR\subseteq bR$ and $cay=c$. 	If $a$ is left and right $(b,c)$ invertible then $a$ is $(b,c)$ invertible. \\
\end{defn}
\noindent When $a$ has a $(b,c)$ inverse (say $y$) then  $bR \cong cR$ (and $Rb \cong Rc$). Here $b,c$ is automatically  von Neumann regular since we can write $b = (y)ab = (buy)ab$ and $c = ca(y)= ca(yvc)$ for some $u,v \in R$. 
\begin{remark}{\label{exchange-bcfiniteness}}
One-sided $(b,c)$-invertibility, unlike $(b,c)$-invertibility, does not force $b$ or $c$ to be von Neumann regular. Rather, it amounts to a comaximality condition: \emph{ $a$ is right $(b,c)$-invertible if and only if $(ab)R+\operatorname{ann}_r(c)=R$, and left $(b,c)$-invertible if and only if $R(ca)+\operatorname{ann}_l(b)=R$}. Now in  \cite[Lemma ~ 20]{siddique1}, we showed that if $R$ is exchange we have

$$\{t\in R : sR+\operatorname{ann}_r(t)=R\} \;=\; \bigcup_{y\in\operatorname{T}(s)} Ry ,$$

\noindent and dually, $Rs+\operatorname{ann}_l(t)=R$ holds  when $t\in yR$ for some $y\in\operatorname{T}(s),$ where $\operatorname{T}(s)$ denotes the set of all outer inverses of $s$. \\
Hence, given  $a \in R$ we can determine from below what kind of values $b,c \in R$ can be chosen to force $a$ to be $(b,c)$ invertible.\\ \emph{ If $R$ is an exchange ring,} then \\
$(i)$ \emph{ $a$ is right $(b,c)$-invertible if and only if $c\in Ry$ for some $y\in\operatorname{T}(ab)$, and in this case $(by)$ is a right $(b,c)$-inverse of $a$. \\
$(ii)$ $a$ is left $(b,c)$-invertible if and only if $b\in yR$ for some $y\in\operatorname{T}(ca)$, in which case $(yc)$ is a left $(b,c)$-inverse of $a$}.\\

\end{remark}

 \noindent In the two-sided case not only are $b,c$ von Neumann regular but also  $ab,ca $ as can be easily seen.  We can then replace $bR$ and $Rc$ by $eR$ and $Rf$ with $e,f \in \operatorname{Idem}(R)$. We assume it in the following result in our block invertibility characterization of $(b,c)$ invertibility.


\begin{prop}{\label{drazin-inverse}}
	Let $a,b,c\in R$ with $b,c$ regular. Then there exists  $e,f\in\operatorname{Idem}(R)$ satisfying $eR=bR$ and $Rf=Rc$. Then
	\begin{enumerate}[label = (\roman*)]
		\item  $a$ is  $(b,c)$-invertible if and only if $f \overset{\scriptscriptstyle fae}{\sim} e$
		\item  $a$ is left $(b,c)$-invertible if and only if $e \lesssim^{\scriptscriptstyle fae} f$;
		\item $a$ is right $(b,c)$-invertible if and only if $f \;{}^{\scriptscriptstyle fae}{\lesssim}\;e $.
	\end{enumerate}
\end{prop}
\begin{proof}
	$(i)$ $(\implies)$: Let $ y$  be the  \emph{Bott-Duffin} $(b,c)$ inverses of $a.$ Since $eR=bR$ and $Rf=Rc$, set $e = br$ and $f = sc$ for $r,s \in R$. Then clearly, $yab = b \implies  yae = e$ and $c = cay \implies f = fay$. Also $y\in bRy\cap yRc \implies y = ey = yf \implies y \in eRf$. \\
	Set $\theta = fae$ and $\theta' = eyf$. Then $f \overset{\scriptscriptstyle fae}{\sim} e$ since  \[ \theta \theta' = (fae)(eyf)= fa(y) = f \,\, \text{ and } \,\, \theta' \theta = (eyf)(fae) =yae= e \]
	
	\noindent $(\Longleftarrow):$ Let $f \overset{\scriptscriptstyle fae}{\sim} e.$ Then there exists $y = eyf \in eRf$ such that $(fae)y = f$ and $y(fae) = e$. Then using $y(fae)y = y$ and  $eR = bR$ and $Rf = Rc$ gives  $y\in bRy\cap yRc.$ Set $b = er$ and $c = wf$. Then clearly $b = (e)r = (yae)r= yab$. Similarly, $c = wf = cay$.\\ Hence $y$ is a $(b,c)$ inverse of $a$.\\
	$(ii)$ $(\implies)$:    Let $ y$  be the left $(b,c)$ inverses of $a.$ Then $Ry \subseteq Rc = Rf \implies yf = y$.  Also, $yab = b$ and $bR = eR$ gives $yae = e$. \\
	Take $\theta = fae$ and $\theta' =  ey$. Then $e \lesssim^{\scriptscriptstyle fae} f$ \[ \theta' \theta = (ey)(fae)= eyae = e \]
	$(\Longleftarrow)$: Suppose $  e \lesssim^{\scriptscriptstyle fae} f$. Then there exists $z \in eRf$ such that $z(fae) = e$. Using $bR= eR$ we get  $ zfae = e \implies zab = b$.  Also $zf = z \implies z \in Rc.$  Thus  $a$ has a  left $(b,c)$-invertible $z$.\\
	$(iii)$ $(\iff)$: Similar to $(ii)$.
\end{proof}

\begin{remark}{\label{aya}}
	In the case when $a$ is $(b,c)$ invertible, the inverse $y$ is also von Neumann regular and is an outer inverse of $a$  since $y \in bRy$ and $yab = b \implies yay = y$. Furthermore, the idempotents $ya$ and $ay$ generate the right and left principal ideals $bR$ and $Rc$ respectively. This is because $b = yab \in (ya)R$ and $ya  \in bR(ya) \subseteq bR$. Same for $Rc$.
\end{remark}
\noindent Thus we have the following proposition. 
\begin{prop}
	Let $R$ be a ring and fix arbitrary $b,c \in R.$ Let $a \in R$. Then the following are equivalent:
	\begin{enumerate}[label = (\roman*)]{\label{bc-invertible}}
		\item $a$ is $(b,c)$-invertible
		\item  There exists $e,f\in \operatorname{Idem}(R)$ with $eR=bR$ and $Rf=Rc$ such that $f\overset{fae}{\sim}e$
		\item There exists $y \in T(a)$  such that  $(ya)R=bR$ and $R(ay)=Rc$.
	\end{enumerate}
	In this case the $(b,c)$ inverse of $a$ is indeed  the inverse of $(fae)$ denoted as $(fae)' \in eRf$.
\end{prop}
\begin{proof}
	Straightforward from Remark \ref{aya} and Proposition \ref{drazin-inverse}. For $(iii)$ note that $ay \overset{aya}{\sim} ya$ always holds for any $y \in \operatorname{T}(a)$.
\end{proof}
Many other well-known inverses of Von Neumann regular elements are particular instances of Drazin's $(b,c)$ inverses  \cite{drazin1}. We note down the block-invertible equivalency of some of those inverses. Here $\operatorname{I(a)}$ denotes the set of all inner inverses of $a$.
\begin{enumerate}
	\item \emph{Group inverse}: The group inverse, if it exists, is the unique $y \in \operatorname{I(a)}$ with $aya=a$, $yay=y$ and $ay=ya$. It is denoted by $a^{\#}$ and is the $(a,a)$-inverse of $a$. \\
	Using Prop. \ref{bc-invertible} $(iii)$, $a^{\#}$ exists if and only if $yaR = aR$ and $Ray = Ra$ for some $y \in \operatorname{T(a)}. $ 
	Since for every $ y \in \operatorname{I(a)}$ we have $(ay)R = aR$ and $R(ya) = Ra, $ then by Prop. \ref{bc-invertible} $(ii)$ \\
	
	 $a^{\#}$ exists if and only if $ya\overset{(ya)a(ay)}{\sim}ay$ for some $y\in I(a)$, and then $a^{\#}=((ya)a(ay))'$.\\

	\item \emph{Moore-Penrose Inverse}: In a ring with involution $*$, the Moore-Penrose inverse of $a$ is the unique $x$ with $axa=a$, $xax=x$, $(ax)^*=ax$ and $(xa)^*=xa$. It is denoted as $a^{\dagger}$  and is the $(a^*,a^*)$-inverse of $a$. Using Prop. \ref{bc-invertible} $(ii)$ \\
	
	 $a^{\dagger}$  exists if and only if $(ay)^*\overset{(ay)^*a(ya)^*}{\sim}(ya)^*$ for some $y\in I(a)$, and  $a^{\dagger}=\bigl((ay)^*a(ya)^*\bigr)'$.\\
\end{enumerate}

Proposition \ref{drazin-inverse} leaves open the question under what circumstances does when one-sided $(b,c)$-invertibility forces two-sided $(b,c)$-invertibility. Answering it requires a notion of stable range for skew corners, which we introduce in Section ~ 3  and return to the question there. Here we examine in detail a special case of the $(b,c)$-inverse, Drazin’s Bott–Duffin inverse \cite{drazin1}. Recall that for idempotents $e,f \in R$, an element $a \in R$ has a unique \emph{Bott-Duffin} $(e,f)$ inverse if there exists a $y\in R$ such that \begin{enumerate}
	\item $y = ey = yf $
	\item $yae = e$ and $fay = f$.\\
\end{enumerate}

From the definition, it can be seen that  for any $a,b,c \in R$, and any pair of idempotents $e,f \in \operatorname{Idem}(R)$ such that  $bR = eR$ and $Rf = Rc$,  the unique Drazin $(b,c)$ inverse (if it exists) co-incides with the Bott-Duffin $(e,f)$ inverse. Therefore we can write for  $e, f \in R$ and $a \in R$: \\
\emph{$a$ has a Bott-Duffin $(e,f)$ inverse if and only if $f \overset{\scriptscriptstyle fae}{\sim} e$}\\


\noindent  In addition if  $(1-e) \sim (1-f)$ then the following holds. 

\medskip 

\begin{theorem}{\label{bott-clean}}
	Let $a \in R$ and idempotents $e^2 = e, f^2 = f$ and  let $\beta$ be an invertible block  in $(1-e)R(1-f)$ with inverse $\beta'$. Then the following are equivalent.
\begin{enumerate}[label = (\roman*)]
	\item $a$ has a Bott-Duffin $(e,f)$ inverse y.
\item $fae \pm \beta'$ is invertible in $R$ (equivalently $ fa \pm \beta'$, or  $ae  \pm \beta'$)
\end{enumerate}
 In this case, $(fae \pm \beta')^{-1} = eyf \pm \beta$.
\end{theorem}
\begin{proof}
$(i) \implies (ii):$ Since $a$ has a Bott-Duffin $(e,f)$ inverse (say) $y$, we have $f \overset{\scriptscriptstyle fae}{\sim} e$ with  $(fae)'=eyf$. Consider the associate pair $(fae + \beta'),(eyf + \beta)$. Then 
\[(fae + \beta')(eyf + \beta) = \underbrace{(fae)(eyf)}_{f} + (fa\underbrace{e)(\beta}_{0}) + \underbrace{(\beta')(e}_{0}yf) + \underbrace{\beta' \beta}_{1-f} = 1\]
Similarly, $(eyf + \beta)(fae + \beta') = e + (1-e) = 1$. \\
Also, $(fae - \beta')(eyf - \beta) = \underbrace{(fae)(eyf)}_{f} - (fa\underbrace{e)(\beta}_{0}) - \underbrace{(\beta')(e}_{0}yf) + \underbrace{\beta' \beta}_{1-f} = 1,$ \\ 
and  similarly $(eyf - \beta)(fae - \beta') =  1$. \\
Thus $fae \pm \beta'$ and $eyf \pm \beta$ are invertible in $R$. \\
\noindent $(ii) \implies (i):$ Suppose $(fae + \beta') \in U(R)$. Set $s = (fae + \beta')^{-1}$ and consider $y = e(s)f$. Then, \[(fae + \beta')(s) = 1 \implies f(fae + \beta')s = f \implies (fae)s = f \implies (fae)y = f\]
and 
\[s(fae + \beta') = 1 \implies s(fae + \beta')e = e \implies s(fae) = e \implies y(fae) = e\]
Thus, $f \overset{\scriptscriptstyle fae}{\sim} e$  and consequently, from above, we get $a$ is \emph{Bott-Duffin} $(e,f)$ invertible.\\

\begin{claim}
  $fae \pm \beta'$ is invertible in $R \iff fa \pm \beta'$ is invertible in $R$.
\end{claim}

Set $\eta = y fa(1-e) \in  eR(1-e)$. Then $\eta^2 = 0 \implies (1+ \eta) \in U(R)$. Then using $\beta' \in (1-f)R(1-e)$ we get
\begin{align*}
	(fae + \beta')( 1 + \eta ) &= fae + fae \eta + \beta' + \underbrace{\beta' \eta}_{0}  = fae + \underbrace{fae(y}_{f}fa(1-e)) + \beta' \\
	&= fae + fa(1-e) + \beta' = fa + \beta'
\end{align*}
Therefore, $fae +\beta' \in U(R) \iff fa + \beta' \in U(R)$. \\
Similarly, $fae -\beta' \in U(R) \iff fa - \beta' \in U(R),$ using the same inverse $( 1 + \eta )$.  

\begin{claim}
  $fae \pm \beta'$ is invertible in $R \iff ae \pm \beta'$ is invertible in $R$.
\end{claim}

	Taking $\gamma = (1-f)aey \in (1-f)Rf$ and continuing similarly as above we get 
	\begin{align*}
	(1+\gamma)(fae \pm \beta') &= fae + \gamma (fae) \pm \beta' \pm \underbrace{ \gamma \beta' }_{0} = fae  + (1-f)a \underbrace{ey (fae)}_{e} \,\,  \pm \,\,  \beta' \\
	&= ae  \pm  \beta'
	\end{align*}
Hence our claim holds true. \\


\noindent The added consequence $(fae \pm \beta')^{-1} = eyf + \beta$  is straightforward from the proof.
\end{proof}

\noindent As an immediate consequence, setting $ e = f$ generates and strengthens the following result from Drazin \cite{drazin1} .\\
\begin{prop}
For any associative ring $R$ with 1 and any $a, e, y \in  R$ with $e^2 = e$, if $a$ has a Bott-Duffin $(e,e)$ inverse,  
 then $1-e + ae$ is invertible, with $(1-e + ae)^{-1} = 1-ay + y$ and $y= e(1-e + ae)^{-1}$. \\
 Dually, also $1-e + ea$ is invertible, with $(1-e + ea)^{-1}=1-ya+y,$ and $y=(1-e + ea)^{-1}e.$
\end{prop}
\begin{proof}
	Taking $\beta = \beta' = (1-e) \in U\bigl((1-e)R(1-e)\bigr)$ in Theorem \ref{bott-clean}, we obtain $ea \pm (1-e) \,\, ,\,\, ae \pm (1-e) \in U(R).$
	Additonally, if $a$ has a Bott-Duffin $(e,e)$ inverse, say $y$ in the skew corner $eRe,$ then we have $y = ey = ye; yae = e = eay$. Then from Theorem \ref{bott-clean} Claim 1, for $\eta = yea(1-e)= ya - e$
\[	\bigl(ea \pm (1-e)\bigr) = (eae \pm (1-e)) (1 + \eta)\]
\begin{align*}
	\implies 	\bigl(ea \pm (1-e)\bigr)^{-1} &= (1-\eta)(eae \pm (1-e))^{-1} \\
	&= (1-\eta)(eye \pm (1-e)) \\
	&=eye  \pm (1-e)  \mp \eta \\
	&= y \pm (1-e) \mp (ya - e) =y \pm (1-ya)
\end{align*}
Similarly, using Claim 2,  $\bigl( ae \pm (1-e) \bigr)^{-1} = y \pm (1-ay)$
\end{proof}

 	Applying left or right $(b,c)$ invertibility on idempotents $e,f$ where $bR = eR$ and $Rc = Rf$  allows us to  define the concepts of right or left Bott-Duffin inverse of an element $a$. 
 
\begin{defn}
	Let $e,f\in R$ be idempotents and let $a\in R$. An element $y\in eRf$  is called
	\begin{enumerate}[label = (\roman*)]
		\item a \emph{right Bott--Duffin $(e,f)$-inverse} of $a$ if $\;fay=f$ so that  $f \; \, \lesssim^{\scriptscriptstyle y }  \; \, e$ 
		\item a \emph{left Bott--Duffin $(e,f)$-inverse} of $a$ if $\;yae=e$ so that  $e \; \; {}^{\scriptscriptstyle y }{\lesssim} \; \; f$
	\end{enumerate}
	 \end{defn}
	If $y_1$ is a left, and $y_2$ a right Bott--Duffin $(e,f)$-inverse	of $a$ resp., then $y_1=y_1f=y_1(fay_2)=(y_1ae)y_2=ey_2=y_2$ and this common element is the Bott--Duffin $(e,f)$ inverse of $a$.   Moreover, if $y\in eRf$ is either a right or left Bott–Duffin $(e,f)$-inverse of $a$, then $y$ is an inner inverse of $fae$ and $a$ is an inner inverse of $y$; that is, $(fae)\,y\,(fae)=fae$ and $yay=y$. This connects one-sided Bott–Duffin invertibility with von Neumann regularity as we  see below. 

 \begin{prop}{\label{bott-regular}}
 Let $R$ be a ring and $a\in R$. Then the following are equivalent:
 \begin{enumerate}[label = (\roman*)]
 	\item $a$ is von Neumann regular 
 	\item There exists idempotents $e,f \in R$ such that $ a \in fRe$ and $a$ has a Bott-Duffin $(e,f)$ inverse.
 \item There exist idempotents $e,f\in R$ with $a\in fRe$ such that $a$ has
 a right Bott--Duffin inverse.
 \item There exist idempotents $e,f\in R$ with $a\in fRe$ such that $a$ has
 a left Bott--Duffin inverse.
 \end{enumerate}
 \end{prop}
 \begin{proof}
$(iii) \implies (i) :$ Suppose there exists idempotents $e, f \in R$ such that  $ a= fae \in fRe$ and $a$ has a right Bott-Duffin $(e,f)$ inverse $y$. Then by definition, $y = ey = yf , fay = f$. Then we have 
 \[  f = fay = (fae)(y) = ay  \text { and }  f(ae) = fae = a  \implies aR = fR\]
Thus $a$ is von Neumann regular.\\
$(i) \implies (ii)$ Conversely, let $a$ be von Neumann regular. Then there exists $x\in R$  such that $a = axa$. Take $y = xax, e = xa$ and $f = ax$. Then $a = (ax)(a)(xa) = fae \in fRe$ and  $ xax = xa (x) ax = exf \in eRf$ and we have $ax \overset{a}{\sim} xa$.
Thus $a \in fRe$ and $a$ has a Bott-Duffin $(e,f)$ inverse $(a')= xax$.\\
$(ii) \implies (iii):$ Straightforward.\\
$(iv) \implies (i):$  Similarly as above. 
 \end{proof}
 \begin{remark}
 We see that for idempotents $e,f \in R$ if  an element $a$ is anchored in the skew-corner $fRe,$ \emph{i.e.} $a = fae,$  then every left (resp. right) Bott-Duffin $(e,f)$ inverse of $a$ is  Bott-Duffin invertible  for some $(e,f)$.
 \end{remark}
In the above characterization, observe that  if there exists idempotents $e,f \in R$ such that $ a \in fRe$ and $a$ has a Bott-Duffin $(e,f)$ inverse, say $y \in eRf$, then \[ay \equiv 1   \pmod{R(1-f)} \text{  and  } ya \equiv 1 \pmod{(1-e)R}\] If every left unit lifts modulo every left ideal, which is a characterization of stable range one in \cite{siddique0} forces $a$ to be unit-regular. Indeed, for $u\in U(R)$ we have $euf=y$ if and only if $e(u-y)f=0$, and since $\forall y \in R, $ $y-eyf = ey(1-f)+(1-e)y, $ therefore  
we get the following characterization for unit-regularity in terms of Bott-Duffin $(e,f)$ inverse.

\begin{prop}\label{bott-unit-regular}
	Let $R$ be a ring and $a\in R$. Then $a$ is unit-regular if and only if
	there exist idempotents $e,f\in R$ and a unit $u\in U(R)$ such that
	$a\in fRe$ and $euf$ is a Bott--Duffin $(e,f)$-inverse of $a$.
\end{prop}

The characterization of unit regularity obtained above is an element-wise statement: it detects, block by block, when a regular element admits a unit inner inverse. The classical theory of unit regularity, going back to Ehrlich and Handelman, is  phrased at the level of idempotents, in terms of the equivalence $e\sim f$. We showed earlier that  equivalence of idempotents is itself a statement about block inverses: $e \overset{x}{\sim} f$ says precisely when the skew corner $eRf$ contains an
invertible block $x$. Using Block generalized inverses, we provide a short proof to Ehrlich \cite{ehrlich} and Handelman's\cite{handelman} classical result on equivalence of idempotents. 
\begin{theorem}[\cite{ehrlich},\cite{handelman}]{\label{ehrlich}}
	\emph{Every von Neumann regular element in a ring R is unit-regular if and only if for every two idempotents $e,f \in R$ , $e \sim f \implies (1-e) \sim (1-f)$.}\\
\end{theorem}
\begin{proof}
	$\Longleftarrow :$ Let $a \in R$ be regular and set $x \in \operatorname{Ref}(a)$. Then choosing $f = ax$ and $e = xa,$ from Prop. \ref{bott-regular}  we obtain $x$ to be the Bott-Duffin $(e,f)$-inverse of $a$ and $f \overset{fae}{\sim } e$. Now $f \sim e \implies (1-f) \sim (1-e)$. Then there exists  invertible block $\beta'$ in $(1-f)R(1-e)$ with inverse $\beta \in (1-e)R(1-f)$, and subsequently by Theorem \ref{bott-clean},  $fae + \beta' = a+ \beta' \in U(R)$ with inverse $x + \beta$. Set $v = (x+\beta)$. Then \[ava = a(x + \beta)a = axa + a \beta a = a + (fae)\beta a = a\]
	
\noindent 	$ \implies :$ Suppose $ f \overset{b}{\sim} e $ for  idempotents $ e,f \in R$ with $b' = x$.  We prove that $(1-f) \sim (1-e)$. \\
	Then clearly $bx = f, xb = e$, and $bxb = b$ and  therefore by assumption, there exists $ u \in U(R)$ such that $b = bub$.  Label $f_1 = bu$ and $e_1 = ub$. Then both $e_1, f_1 \in \operatorname{Idem}(R)$.\\
	Consider the elements $\gamma = (1-e_1)u(1-f_1)$ and $\gamma' = (1-f_1)u^{-1}(1-e_1)$. Then it can be seen that 
	\[(1-e) \overset{(1-e_1)}{\sim} (1-e_1) \overset{\gamma}{\sim }  (1-f_1) \overset{1-f_1}{\sim} (1-f)\]

\end{proof}

Now, in the proof of Proposition~\ref{bott-regular} on a characterization of von Neumann regularity, we observed when  a subequivalence $\lesssim$   improves to an equivalence $\sim$. It is therefore natural to ask what becomes of Theorem~\ref{ehrlich} when the hypothesis $e\sim f$ is weakened to $e\lesssim f$. \\
Observe first that the conclusion must reverse direction: That is, if for idempotents $e,f \in R$ we have $e\leq f$  then automatically, $1-f\leq 1-e.$  Therefore our complement condition compatible with $\lesssim$ is $1-f\lesssim 1-e$. Observe that such a ring $R$ where for arbitrary $e,f \in R,$  $e \lesssim f \implies (1-f) \lesssim (1-e)$, is always directly finite. This is because  whenever $pq = 1$ for some $p,q \in R$ then $1 \overset{p}{\sim} qp$ which gives $(1-qp) \lesssim 0$ forcing $qp = 1$. \\
Surprisingly, the weakened hypothesis loses none of its strength as we see below.
\begin{theorem}
	\emph{Every von Neumann regular element in a ring R is unit-regular if and only
		if for arbitrary   idempotents $e,f \in R,$  $e \lesssim f \implies (1-f) \lesssim (1-e)$.}\\
\end{theorem}
\begin{proof}
	$ \implies :$ Let every von Neumann regular element in $R$ be unit-regular. Choose arbitrary idempotents $e,f \in \operatorname{Idem}(R)$ such that $e \lesssim f$. \\
	Then by definition, there exists an idempotent $ p \in \operatorname{Idem}(R)$ such that $e \sim p$ and $ p \leq f$. Now by Theorem \ref{ehrlich}, $p \sim e \implies (1-p) \sim (1-e)$. 

	\noindent Also, since $ p \leq f $ then by definition, $ p = fp = pf$. Therefore,  $(1-p)(1-f) = 1-p-f + pf = (1-f)$ and $(1-f)(1-p) = 1-p-f + fp = (1-f) \implies (1-f) \leq^{-} (1-p)$.  Then $ (1-f)  \lesssim (1-p)$. \\
	Thus by Lemma \ref{chotabhim}$(1)$,  $(1-f)\lesssim (1-p) \sim (1-e) \implies (1-f) \lesssim (1-e) $.\\
	\noindent Conversely, let $a \in R$ be regular and set $x \in \operatorname{Ref}(a)$. Then from Prop. \ref{bott-regular} $(iii)$  choosing $f = ax$ and $e = xa$ we obtain $f \; \; {}^{fae}\lesssim \; e$ where $a = fae$ and $x \in eRf$. \\
	Then by hypothesis,  $(1-e) \;\; {}^{r}\lesssim \;(1-f)$ where $r \in (1-e)R(1-f)$ with $r' \in (1-f)R(1-e)$ such that $rr' = 1-e$. Then $(x+r)$ is right-invertible in $R$ since $(x+r)(a+r')= e+(1-e) = 1$ and therefore using direct finiteness as shown earlier, $x+r \in \operatorname{U}(R)$. Label $u = x+r$. Then $aua = a(x+r)a= axa =a. $
\end{proof}

	In both global characterizations of unit-regularity above, the implication  required is for every pair of idempotents of $R.$ However it is of interest to know when such an equivalence holds true locally for fixed arbitrary idempotents $e,f$. The relations $e\sim f$ and $(1-e)\sim(1-f)$ say precisely that the complementary skew corners $fRe$ and $(1-f)R(1-e)$ contain invertible blocks. In the next couple of results we approach the question through (one-sided) block inverses in  their respective skew specific corners.\\

We begin with what the two relations give when both hold. The Lemma below allows us to manufactures units in situations where a direct computation becomes very technical. The mechanism is the ring-theoretic counterpart of  block inversion in matrices. Compressing $x\in R$ on the left by idempotents $p,1-p$ and on the right by idempotents $e,1-e$ produces the  \emph{two-sided Peirce decomposition}

\[x \;=\; \underbrace{pxe}_{\alpha} \;+\; \underbrace{px(1-e)}_{\eta}
\;+\; \underbrace{(1-p)xe}_{\gamma} \;+\; \underbrace{(1-p)x(1-e)}_{\beta},
\]
which we can represent in matrix form
\[
x \equiv \begin{pmatrix} \alpha & \eta \\ \gamma & \beta \end{pmatrix}
\]
The two diagonal entries lie in \emph{different} skew corners, so \emph{invertible} is  read as an invertible block in the corresponding skew corner. This is the degenerate case of Banachiewicz's \cite{bana} full inversion formula for a partitioned matrix where either of the off-diagonal elements $B,C$ are zero. In a matrix $M=\left(\begin{smallmatrix} A & C\\ B & D\end{smallmatrix}\right)$, if $A$ is invertible and either of $B$ or $C$ is zero, then
$M$ is invertible if and only if the Schur complement $S=D-BA^{-1}C = D$ is, in which
case
\[
M^{-1}=\begin{pmatrix}
	A^{-1} & -A^{-1}CS^{-1}\\
	-S^{-1}BA^{-1} & D^{-1}
\end{pmatrix}.
\]
where either $-A^{-1}CS^{-1} = 0$ or $-S^{-1}BA^{-1} = 0,$ and  where $A$ and $D$ are invertible matrices of different  sizes. Thus, in Peirce form, 
if $\alpha$ is an invertible block of $pRe$ with inverse $\alpha'\in eRp$,
and the Schur complement $S=\beta\in(1-p)R(1-e)$ is an invertible block
with inverse $\beta'\in(1-e)R(1-p)$, then $x\in U(R)$ whenever
$\gamma=0$ or $\eta=0$, and
\[
x^{-1}=\alpha'-\alpha'\eta\beta'-\beta'\gamma\alpha'+\beta'.
\]

 \noindent We shall refer to the following lemma as the  \emph{ Block Lemma} throughout this article. 
\begin{lemma}\label{triangular}
	
 Let $p,e$ be arbitrary idempotents and $p \overset{\alpha}{\sim} e$ and $(1-p) \overset{\beta}{\sim}  (1-e)$. Then the following results hold true. 
	\begin{enumerate}[label = (\roman*)]
		\item $\alpha+\gamma+\beta  \in U(R)$ for every $\gamma \in (1-p)Re.$ 
		\item $\alpha+\eta+\beta  \in U(R)$ for every $\eta \in pR(1-e).$ \\
	\end{enumerate}

	
\end{lemma}
\begin{proof}
	$(i) \quad$		Let $(\alpha, \alpha')$ be an invertible block in $(pRe,eRp)$ and $(\beta, \beta')$ be an invertible block in \\
	$(1-p)R(1-e), (1-e)R(1-p)$. Then for any $\gamma \in (1-p)Re$  the cross products vanish 
	\[	\alpha\beta'=0=\beta\alpha' ; \alpha'\beta=0=\beta'\alpha ; \alpha'\gamma=0 = \gamma\beta'\]
	Set $d = (\alpha + \gamma + \beta)$. Then,	since $\gamma=(1-p)\gamma$ we have 
	\begin{align*}
		d(\alpha'+\beta'-\beta'\gamma\alpha') &= \alpha\alpha'+\gamma\alpha' + \beta \beta'- \underbrace{\alpha \beta'}_{0}\gamma \alpha' - \underbrace{\gamma \beta'}_{0}\gamma \alpha'-  \beta\beta'\gamma\alpha' \\
		&= p + \gamma\alpha' + (1-p) - (1-p)\gamma \alpha' \\
		&= p + (1-p) + \left[\gamma \alpha' - \gamma \alpha'\right] = 1
	\end{align*}

	\noindent On the other hand, since $\gamma=\gamma e$ we also get 
	\begin{align*}
		(\alpha'+\beta'-\beta'\gamma\alpha')d &= \alpha'\alpha+ \beta'\gamma + \beta' \beta- \beta' \gamma (\alpha' \alpha) - \beta' \gamma\underbrace{\alpha' \gamma}_{0}\gamma \alpha' - \beta' \gamma\underbrace{\alpha'\beta}_{0} \\
		&= e + \beta'\gamma + (1-e) - \beta' \gamma(e) = 1 
	\end{align*}

\noindent 	Thus, $d = \alpha+\gamma+\beta  \in U(R)$ for every $\gamma \in (1-p)Re$. \\

\noindent 	$(ii) \quad$ The proof of $(ii)$ follows similarly as above.
\end{proof}

As an immediate consequence for $\gamma = 0$, if $ p \overset{\alpha}{\sim} e$ and $ (1-p) \overset{\beta}{\sim} (1-e)$ 	then $\alpha+\beta\in U$ with $(\alpha+\beta)^{-1}=\alpha'+\beta'$. 	In the diagonal case where $p=e$, we have  $\alpha\in pRp$ and $\beta\in (1-p)R(1-p), $ then 	$\alpha + \beta\in U$ if and only if $\alpha\in U(pRp)$ and $\beta\in U\left((1-p)R(1-p) \right)$.\\
Note that Blocks Lemma \ref{triangular} above is a generalization of Lemma 6.1 in \cite{khu2} Khurana and Lam. 

\begin{lemma}[\cite{khu2} Khurana and Lam]
	Let $e, f \in  idem (R) $ be such that $e+f = 1$, and let $x \in U(eRe)$, $y \in U(f Rf )$.
	Then $x + p + y \in U(R)$ for every $p \in f R e$.
\end{lemma}
\bigskip

\noindent  It is natural to ask to what extent the converse holds: if such an $x$ is a unit, must the diagonal blocks $\alpha$ and $\beta$ themselves admit block inverses? Observe that for some distinct $p,e \in \operatorname{Idem}(R),$ whenever $x=\alpha+\eta+\beta \in U(R)$ where $\alpha\in pRe, \eta \in pR(1-e), $  and  $\beta\in(1-p)R(1-e)$, then $\alpha\in pRe$ is  automatically left block-invertible and $\beta$ is right block invertible in their respective skew-corners since $(ex^{-1}p)\alpha = e$ and $\beta \left((1-e)x^{-1}(1-p)\right) = 1-p$ .\\ In general the blocks $\alpha$ or $\beta $ is are not both-sided invertible in their respective corners, and  the unit condition alone does not force the invertibility of either of those blocks. What it does force is an association between the two $: \alpha$ now admits a right block inverse precisely when $\beta$ is left-block invertible, and vice versa. In other words $p \overset{\alpha}{\sim} e   \iff (1-p) \overset{\beta}{\sim} (1-e)$. The following lemma makes this relationship precise.  

\medskip 
\begin{prop}{\label{converse}}
Let $\alpha\in pRe$ and  $\beta\in(1-p)R(1-e)$. Then the following  hold. 
\begin{enumerate}
	\item If for some $\eta \in pR(1-e), $ we get $\alpha+\eta+\beta\in U(R)$ then 

 \[ p \;  \mathrel{{}^{ \alpha}\!\lesssim} \; e \iff (1-e) \lesssim^{\beta} (1-p) \]

\item If for some $\gamma \in (1-p)Re, $ we get $\alpha+\gamma+\beta\in U(R)$ then 

 \[ e \;  \lesssim^{\alpha}\; p \iff (1-p)  \mathrel{{}^{ \beta}\!\lesssim}  (1-e) \]

\end{enumerate}
\begin{proof}
	$(1)$	Label $x = \alpha+\eta+\beta \in U(R)$ for some $\eta \in pR(1-e).$\\
	$(\implies) :$ Let $p \;  \mathrel{{}^{ \alpha}\!\lesssim} \; e$ with $\alpha' \in eRp.$  Then $\alpha \alpha' = p .$ We prove that $\beta’ =(1-e)x^{-1}(1-p)$ satisfies $\beta\beta’=1-p$ and $\beta’\beta=1-e$. 
	\begin{align*}
	\text{Now, } x\alpha' &= \left(\alpha+\eta+\beta\right) \alpha' \\
	& = \alpha\alpha' + \eta\alpha' + \beta\alpha' \\
	&= p + \bigl(\eta(1-e)\bigr)(e\alpha’) + \bigl(\beta(1-e)\bigr)(e\alpha’)=p \\
	\implies \hspace{3mm} \alpha' &=  x^{-1}p
	\end{align*}
	Now, $e\alpha' = \alpha' \implies (1-e) \alpha' = 0 \implies (1-e)x^{-1}p = 0$.\\
	
\noindent 	Since $\eta = p\eta$ and $\beta = (1-p)\beta$ we get $x (1-e) = \left(\alpha+\eta+\beta\right) (1-e) = \alpha(1-e) + \eta(1-e)+ \beta(1-e)=  \eta +\beta $
		\begin{align*}
	\implies \hspace{3mm} (1-e)x^{-1}\left(x (1-e) \right) &=  (1-e)x^{-1} \left(\eta + \beta\right) \\
	&=  (1-e)x^{-1} \left(p\eta + (1-p)\beta\right) \\
	&= \underbrace{(1-e)x^{-1} p}_{0}\eta  + \underbrace{(1-e)x^{-1} (1-p)}_{\beta'}\beta = \beta' \beta
		\end{align*}
		Thus, we have $(1-e) \lesssim^{\beta} (1-p)$ with $\beta'\beta = 1-e$. \\
In fact, since $\beta$ is both left and right invertible,  $\beta$ is block invertible with inverse $\beta' $.
\begin{remark}
\noindent In the proof above, we simply assumed $\alpha$ is right invertible whence $\alpha \alpha' = p.$ Then as a consequence we not only obtain $\beta$ being an invertible block in $(1-p)R(1-e)$ but $\alpha$ being also  left invertible. Indeed, $e = (eM^{-1})(Me) = (eM^{-1}) (\alpha + \eta + \beta)e = (eM^{-1}) \alpha$. Also $M \alpha' = p \implies M\alpha' \alpha = p \alpha = \alpha \implies (\alpha'\alpha) = M^{-1} \alpha.$ Therefore $e = e M^{-1} \alpha = e (\alpha' \alpha) \implies \alpha' \alpha = e$.\\
\end{remark}
\noindent $(\Longleftarrow):$ Suppose $(1-e) \lesssim^{\beta} (1-p)$ with $\beta'\beta = 1-e$ for $\beta' \in (1-e)R(1-p)$.
We prove that $\alpha’ =eM^{-1}p$ satisfies $\alpha\alpha’=p$.\\
Since $\alpha = p\alpha,$ $\, \eta = p \eta $ and $\beta' = \beta' (1-p)$ we have
			\begin{align*}
\beta' M &= \beta' \left(\alpha+\eta+\beta\right)  \\
			& = \beta' \alpha+ \beta' \eta  + \beta'\beta \\
			&=  (\beta' (1-p))(p\alpha) + (\beta'(1-p)p\eta) + (1-e)= (1-e)\\
			\implies \hspace{3mm} \beta' &=  (1-e)M^{-1}
		\end{align*}
		
	Also, $\beta' = \beta'(1-p) \implies \beta' p= 0 \implies (1-e)M^{-1}p = 0$.\\
	
	\noindent 	Now $pM  = (p)\left(\alpha+\eta+\beta\right)  = \alpha + \eta $. Since  $\eta = \eta(1-e)$ and $\alpha = \alpha e$ we  have
			\begin{align*}
		\implies \hspace{3mm} \left(pM \right)M^{-1}p &=   \left(\alpha + \eta\right) M^{-1} p\\
		&=   \left(\alpha(e)  \,\,+ \,\, \eta(1-e)\right)M^{-1}p \\
		&= \alpha\underbrace{eM^{-1} p}_{\alpha'}  \,\,+ \,\,\eta\underbrace{(1-e)M^{-1} p}_{0} = \alpha \alpha'
	\end{align*}
	We therefore have $\alpha \alpha' = p$. \\
	In fact, since $\alpha$ is both left and right invertible,  $\alpha$ is block invertible with inverse $\alpha' $.\\
	
\noindent $(2)$ The proof for $(2)$ holds similarly as above for $\gamma \in (1-p)Re$ by taking $\alpha_1 = \beta, \beta_1 = \alpha, p_1 = 1-p$ and $e_1 = 1-e$ and applying $(1)$.
\end{proof}
\end{prop}
	

\medskip 
\medskip 

 \section{Block Invertibility and Stable Range}
 \bigskip 
 
 We now look at defining the notion of  stable range on skew-corners such as $eRf$ for distinct idempotents $e,f$, since they are intimately related to the concept of invertible blocks.
 
 \begin{defn} Let $e,f \in \operatorname{Idem}(R)$.
The skew-corner $eRf$ has right  stable range one, written $\operatorname{sr}_r(eRf)=1$, if whenever $a\in eRf$ and $b\in eRe$ satisfy $aR+bR=eR$, there exists $y\in eRf$ such that $a+by$ is right block invertible in $eRf$, \emph{i.e.} there exists $z \in fRe$ such that $(a+by)z = e$.\\

 \noindent Similarly, the skew-corner $eRf$ has left stable rank one, written $\operatorname{sr}_l(eRf)=1$, if whenever $a\in eRf$ and $b\in fRf$ satisfy $Ra+Rb=Rf$, there is $y\in eRf$ such that $a+yb$ is left block invertible in $eRf$, \emph{i.e.} there exists $z \in fRe$ such that $z(a+yb) = f$.
 \end{defn}
 \noindent When $e=f$ we obtain the usual notion of stable range on the corner ring $eRe$.
\begin{remark}
Note that $\operatorname{sr}_r(eRf) = 1$ is equivalent to the notion that whenever $ax + b = e$ for $a\in eRf, x \in fRe$ and $b \in eRe$,  there exists $ y \in eRf$ such that $a + by $ is  right block invertible in $eRf.$ This coincides with Ara and Goodearl's \cite{goodearl1} definition of (right) stable rank $1$ on $eRf$. The left stable range condition on $eRf$ can be similarly phrased as follows: \\
 $\operatorname{sr}_l(eRf)=1$ if whenever for $a\in eRf, x \in fRe$ and $b\in fRf$ we have $xa + b = f,$ there exists $y \in eRf$ such that $a + yb$ is left block invertible in $eRf$.
\end{remark}
In general, the right and left stable range conditions on a skew-corner need not be symmetric. This is because if $\operatorname{sr}_r(eRf)=1$ then taking $a = 0$ and $b = e$ gives $y \in eRf$ such that $by$ is a right block invertible in $eRf$, or in other words $e  \;\mathrel{{}^{by}\!\lesssim}  \; f$. But if $\operatorname{sr}_l(eRf)=1$, then taking $a=0$ and $b=f$ gives $y\in eRf$ such that $yb$ is left block invertible in $eRf$, so that $f\lesssim^{yb} e$. We have seen before  with $R=M_2(k)$ with $f=I_2$ and $e=\begin{pmatrix}1&0\\0&0\end{pmatrix} $ that one does not necessarily imply the other. However, we have the following  association. \\
\begin{prop}{\label{leftrightcorner}}
Let $e,f \in \operatorname{Idem}(R)$. Then $eRf$ has right (resp. left) stable range one if and only if its associate skew corner $fRe$ has left (resp. right) stable range one.
\end{prop}
\begin{proof}
It suffices to prove that $\operatorname{sr}_l(fRe) = 1$ when $\operatorname{sr}_r(eRf) = 1$. \\
Suppose $xa + b = e$ for some $a \in fRe, x \in eRf$ and $b \in eRe$. Then there exists $y \in eRf$ such that $x + by$ is right block invertible in $eRf$. Choose $z \in fRe$ such that $(x+by)z = e$.\\
We need a $y' \in fRe$ such that $a + y'b$ is left block invertible in $fRe$.\\
Define $y' = (f-ay)z$ and $s = x + y(f-ax)$. Then $sa = xa+ ya -yaxa =  xa+ yab$ and 
\begin{align*}
s(f-ay) = \left(x + y(f-ax)\right) \left( f-ay\right) &= x - xay + y - yax - yay +yaxay \\
&= \left( x - xay + y \right) - ya \left(x - xay + y  \right) \\
&= \left(e-ya\right) \left(x - xay + y\right)
\end{align*}
where $\left(x - xay + y\right) = x + (e-xa)y = x+ by$. \\
Then $s(a +y'b) = sa + sy'b = sa + (e-ya)\underbrace{(x+by)z}_{e}b = sa + (e-ya)b = xa+b = e$ \\
Thus $a + y'b$ is left block invertible in $fRe$. \\
\end{proof}

\noindent Using  $e = f$, we recover the following result of Vaserstein.
\begin{cor}[Vaserstein~ \cite{vaserstein}]
	A ring $R$ has left stable range one if and only if it has right stable range one. 
\end{cor}

\noindent Right ( resp. Left) stable range one of a skew-corner ring $eRf$ being a stronger form of $e \lesssim f$ (resp. $f \lesssim e$), provides the right conditions under which a subequivalent pair of idempotents gets promoted to equivalency as we see below.  \\
\begin{theorem}
	Let $a \in R$ and $b,c \in R$ be von Neumann regular with $bR = eR$ and $Rc = Rf$ for some $e,f \in \operatorname{Idem}(R)$. Then the following holds. 
	\begin{enumerate}[label = (\roman*)]
		\item If $\operatorname{sr}_l(eRf) = 1$ then $a$ is  left $(b,c)$ invertible if and and only if  $a$ is   $(b,c)$ invertible.
		\item If $\operatorname{sr}_r(eRf) = 1$ then $a$ is  right $(b,c)$ invertible if and and only if  $a$ is   $(b,c)$ invertible.
		\item If $e \sim f$ and $\operatorname{sr}_l(eRe) = 1$  then every one-sided  left $(b,c)$ inverse is two-sided.  
	\end{enumerate}
\end{theorem}
\begin{proof}
\noindent 	$(i)$ Suppose $\operatorname{sr}_l(eRf) = 1$. Under this scenario, we  observe that every right invertible block in $eRf$ is  invertible.\\
	\noindent  Let $a \in eRf$ be right invertible in $eRf$ {\em i.e.} there exists $x \in fRe$ such that $ax= e$. We can write $xa + b = f$ where $b = (f-xa) \in fRf$ and $xa \in fRf$. By definition, there exists $y \in eRf$ such that $a + yb $ is left invertible in $eRf$ {\em i.e.}  there exists $z \in fRe$ such that $z\left(a+ yb\right)= f$. Post-multiplying by $x$ gives $zax = x$ or $z = x$. Thus $x$ is an invertible block in $fRe$ and consequently, $xa = e$. Thus $a$ is an invertible block in $eRf$.\\
	Now if an element $a$ is left $(b,c)$ invertible then by Proposition \ref{drazin-inverse} $(ii)$ we have $e \lesssim^{\scriptscriptstyle fae} f$, and therefore there exists $w \in eRf$ such that $w(fae) = e$. Then $w$ being right invertible in $eRf$ is invertible and therefore, $f \overset{fae}{\sim} e$ which by Proposition \ref{drazin-inverse} $(i)$  proves that $a$ is  $(b,c)$ invertible.\\
\noindent	$(ii)$ If $\operatorname{sr}_r(eRf) = 1$ then we can similarly show that  every left invertible block in $eRf$ is  invertible. And the the proof follows similarly as above. \\
	\noindent $(iii)$  If $e \sim f$ and $\operatorname{sr}_l(eRe) = 1$ then by  \cite[Lemma ~ 2(b)]{goodearl1},  $\operatorname{sr}_r(eRf) = 1$ and $\operatorname{sr}_l(eRf) = 1$. Then using the arguments above we can see that every left and right invertible block in $eRf$ is invertible. Thus every  one-sided  left $(b,c)$ inverse is two-sided.  
\end{proof}
\begin{remark}{\label{direct-finiteness}}
		From $(iii)$ and from \cite[Lemma ~ 2(b)]{goodearl1}, we can see that $\operatorname{sr}_r(eRf) = 1$ and $\operatorname{sr}_l(eRf) = 1$ if and only if $e \sim f$ and $\operatorname{sr}(eRe) = 1$.\\
	\end{remark}
\noindent In the proof above we showed that $\operatorname{sr}_l(eRf) = 1 \left(\text{resp.  } \operatorname{sr}_r(eRf) = 1 \right)  $ then every  right (resp. left) invertible block in $eRf$ is  invertible. Now	recall that if $e \sim f$, every right (resp. left) invertible block of $eRf$ is exactly $x$ times a right (resp. left) block invertible element of $fRf$. The same works for $eRe$ but by post-multiplication by $x$. Hence we  also have a weaker characterization of the form: \\

\noindent 	\emph{If $e \sim f$ and $eRe$ is directly finite then every left  and right  invertible block in $eRf$ is  invertible .  } \\

\noindent The converse does not necessary hold unless $eRf$ contains atleast one left or right inveritble block {\em i.e.}
\emph{If $e \lesssim f$ or $f \lesssim e$, then every left and every right invertible block in $eRf$ is invertible if and only if $e \sim f$ and $eRe$ is directly finite.}\\



\noindent This inference gives  a partial answer to a question \cite[Sec. ~4, Pg 71]{drazin1}: \\

\emph{ If a ring $R$ has stable range one, is it $(b,b)$-finite for all $b \in R$  }? \\

\noindent In \cite{drazin1}, for a fixed $b \in R$, Drazin defines a ring $R$ to be $(b,b)$-\emph{finite}  if for
all $a\in R$, \\$a$ is left $(b, b)$-invertible if and only if $a$ is right $(b, b)$-invertible. \\

\noindent For von Neumann regular elements $b \in R$ the result is true  for directly finite rings.
\begin{prop}
$R$ is directly finite if and only if $R$ is $(b,b)$-finite for every von Neumann regular element $b \in R$.
\end{prop}
\begin{proof}
$(\implies):$ Suppose $R$ is directly finite. Let $b \in R$ be von Neumann regular and $b= bxb$ for some $x\in R$. Set $e = bx$ and $f = xb$. Then $e,f \in \operatorname{Idem}(R)$ such that $bR = eR$ and $Rb = Rf$ and $e \sim f$. 
Then from Prop. \ref{drazin-inverse}$(ii)$, $a$ is left $(b,b)$-invertible if and only if $e\lesssim^{fae}f$, that is, the block $fae\in fRe$ has a left inverse in $eRf$. Likewise, $a$ is right $(b,b)$-invertible if and only if $f\mathrel{{}^{fae}\!\lesssim}e$, that is, $fae$ has a right inverse in $eRf$.
Since $R$ is directly finite, so is its corner ring $eRe$. As $e\sim f$,   Remark \ref{direct-finiteness} shows that any block in $fRe$ with a one-sided inverse in $eRf$ is invertible. Hence, if $a$ is left or right $(b,b)$-invertible, then $fae$ is invertible, that is, $f\overset{fae}{\sim}e$. By Prop. \ref{drazin-inverse} $(i)$, $a$ is then $(b,b)$-invertible, and in particular both left and right $(b,b)$-invertible. Thus $R$ is $(b,b)$-finite.\\
\noindent Conversely, if suppose $ay = 1$ for some $y \in R$. Then $a$ has a right $(1,1)$ inverse $y$. Then $R$ being $(1,1)$ finite forces $ya = 1$. Thus $R$ is directly finite. 
\end{proof}


\noindent In the case when $b$ is not assumed to be regular we obtain  an equivalent characterization of left and right $(b,b)$ invertibility  from  Remark \ref{exchange-bcfiniteness} in Section ~2.\\
If $R$ is an exchange ring\\
\noindent $(i)$ \emph{ $a$ is right $(b,b)$-invertible if and only if $b\in Ry$ for some $y\in\operatorname{T}(ab)$, and in this case $(by)$ is a right $(b,b)$-inverse of $a$. \\
	$(ii)$ $a$ is left $(b,b)$-invertible if and only if $b\in yR$ for some $y\in\operatorname{T}(ba)$, in which case $(yb)$ is a left $(b,b)$-inverse of $a$}.\\

\begin{prop}
	Let $R$ be an exchange ring. Then $R$ is directly finite if and only if $R$ is $(b,b)$-finite for all $b \in R$.
\end{prop}
\begin{proof}
$(\implies):$ Suppose $R$ is directly finite and fix $b \in R$. For any $a \in R,$ let $a$ be right $(b,b)$ invertible. Then from above, $b\in Ry$ for some $y\in\operatorname{T}(ab)$, and in this case $(by)$ is a right $(b,b)$-inverse of $a$. \\
\noindent  Set $w = by$. Then $waw = w$ and $baw = baby = b$. Label $e = wa$ and $f = aw$. Clearly $e,f \in \operatorname{Idem}(R)$ and $f \overset{fae}{\sim} e$. We prove that $bR = eR$ and $Rb = Rf$.\\
Clearly,  $e = bya \implies e  \in bR$ and $ b = bf \implies b \in Rf$.\\
 Now $e =wa = bya \implies e = (ebf)(fyae)$. Label $x = ebf$. Then $x$ is a right-invertible block in $eRf$. Since $eRe$ is directly finite and $e \sim f$, by Remark \ref{direct-finiteness}, $x$ is invertible and subsequently, $(fyae)(ebf)=f \implies fyae\underbrace{baw}_{b} = f \implies f \in Rb$. Therefore $Rb = Rf$.\\
 Also,  $b = bf = b(fyae)(ebf) = \underbrace{bya}_{e}(ebf) \in eR.$ Therefore $eR = bR$.\\
 Hence, by Prop. \ref{drazin-inverse} $(i)$, $a$ is  $(b,b)$-invertible with inverse $w$.\\
 The proof works similarly if $a$ is assumed to be left $(b,b)$ invertible. \\
 Since $a$ is arbitrary,  $R$ is $(b,b)$-finite.\\
 The converse is straightforward from Proposition 3.7.
\end{proof}
\begin{remark}
Without the hypothesis of exchange on $R$, if for some $a\in R$ we choose $b\in R$ such that $ab^2,b^2a$ are von Neumann regular and $R$ is directly finite, then the same technique can be applied alongwith the equality in \cite[Pg~12]{siddique1} to prove that $a$ is left $(b, b)$-invertible if and only if $a$ is right $(b, b)$-invertible. \\

\end{remark}

When $e = f$, the skew corner $eRf$ reduces to the corner ring $eRe$. We now record several cases in which an element of $eRe$ has stable range one. By $	\operatorname{sr}_{\scriptscriptstyle S}(a)$ in a ring $S$ we mean the left stable range, \emph{i.e} whenever $Sa +Sb = S$ for any $b\in S$, there exists  $s\in S$ such that $a + sb \in U(S)$. This is equivalent to the right stable range of the element \cite{khu2}.
We note a simple well-known characterization of  stable range one condition on elements  first proved by Menal and Moncasi in \cite{menal} as pointed out by P. Ara. It turns out to be quite  useful as it is not based on a condition and thus limits the number of independent variables in the pot. Suprisingly, a variation of it also holds for regular and unit-regular elements \cite{lam2}. 

\begin{lemma}{\label{stable}}
	Let $R$ be a ring with identity $1_{\scriptscriptstyle R}.$ Then  the following are equivalent
	\begin{enumerate}[label = (\roman*)]
		\item $ \operatorname{sr}_{\scriptscriptstyle R}(a)=1 $
		\item For  each $s \in R,$ there exists $ r \in R$ such that $a + r(1_{\scriptscriptstyle R}-sa) \in U(R)$ 
		\item For  each $s \in R,$ there exists $ r \in R$ such that $a + (1_{\scriptscriptstyle R}-as)r \in U(R)$ 
	\end{enumerate}
\end{lemma}
\begin{proof}
	$(i) \implies (ii) :$ Let $a$ have stable range one. Choose arbitrary $s \in R$. Then $sa + (1-sa) = 1$. Since $ \operatorname{sr}_{\scriptscriptstyle R}(a)=1, $ we have $a + r(1-sa) \in U(R)$ for some $r\in R$. \\
		$(ii) \implies (i) :$  Conversely, for each $s \in R$ let  $a + r(1-sa) \in U(R)$  for some $r \in R$.  Suppose $xa + b = 1.$ Then $b = (1-xa)$. Take $s = x$. Then $a + r_x(1-sa) = a + r_x(1-xa) = a + r_xb \in U(R)$ for some $r_x\in R$. Thus $a$ has left stable range one.\\
		$(i) \iff (iii):$ Similar, using the left-right symmetry of stable range one elements.
\end{proof}

\medskip 
Khurana and Lam developed the Suspension Theorem in \cite{khu2} that allows transfer of stable range one from the Peirce corner to the ambient ring. We extend this notion from a single Peirce decomposition to transporting stable range across two-sided Peirce decompositions for distinct idempotents. The use of distinct idempotents in the preceding theorem is motivated by  the rectangular-corner techniques of Ara and Goodearl  \cite{goodearl1} . In their study of stable rank of full corners, Ara and Goodearl introduced a stable rank condition for a skew corner $pRq$, and showed that this condition is preserved upon replacing either endpoint idempotent by an equivalent one (\cite{goodearl1}, Lemma 2(b)). \\
 In keeping with the terminology used by Khurana and Lam, we shall refer to the following result as the \emph{General Suspension Theorem}. 
\smallskip 


\begin{theorem}[\textbf{General Suspension Theorem}]{\label{suspension}}
Let $p,e\in\mathrm{idem}(R)$ such that \[p \overset{\theta}{\sim} e \;\; \text{  and  } \;\; (1-p) \overset{\beta}{\sim} (1-e)  \]  Fix arbitrary $a \in eRe, c \in pRp.$ Then the following conditions hold.
\begin{enumerate}
		\item  For every $\eta \in pR(1-e) $ and $\tilde{\eta} \in eR(1-p) $, the following hold.
	\begin{enumerate}[label =(\roman*)]
		\item  $	\operatorname{sr}_{pRp}(c)=1\iff \operatorname{sr}_{\scriptscriptstyle R}\bigl(c\theta +\eta +\beta\bigr)=1 $
		\item  $	\operatorname{sr}_{eRe}(a)=1\iff \operatorname{sr}_{\scriptscriptstyle R}\bigl(a\theta' + \tilde{\eta}  +\beta'\bigr)=1 $
	\end{enumerate}
\medskip
	\item  For every $\eta \in pR(1-e) $ and $\tilde{\eta} \in eR(1-p) $, the following hold.
	\begin{enumerate}[label =(\roman*)]
	\item  $	\operatorname{sr}_{pRp}(c)=1\iff \operatorname{sr}_{\scriptscriptstyle R}\bigl(\theta'c + \tilde{\eta} +\beta'\bigr)=1 $
	\item  $	\operatorname{sr}_{eRe}(a)=1\iff \operatorname{sr}_{\scriptscriptstyle R}\bigl(\theta a + \eta+\beta \bigr)=1 $
\end{enumerate}
\end{enumerate}
\noindent In each of the cases, necessarily $\operatorname{sr}_R(a)=1$.
\end{theorem}

\smallskip
\begin{proof}
	 By definition,   $\theta' \in eRp$ and $\beta' \in (1-e)R(1-p).$\\
$(1) (i)$ For arbitrary $c \in pRp,$ we prove that for every $\eta \in pR(1-e) $ 
\[		\operatorname{sr}_{pRp}(c)=1\iff \operatorname{sr}_{\scriptscriptstyle R}\bigl(c\theta +\eta +\beta\bigr)=1 \]
$\implies :$ Let $	\operatorname{sr}_{pRp}(c)=1$ and label $M = c\theta +\eta +\beta$ for arbitrary $\eta \in pR(1-e)$. Fix an arbitrary $s\in R.$ We prove $M + r(1-sM) \in U(R)$  for some $r \in R$. 
	
	\noindent Consider the element $\tilde{s}=\theta s p \in pRp$. Then in the ring $pRp$ we have  $\tilde{s}c + (p-\tilde{s}c ) = p.$ Hence by $\operatorname{sr}_{pRp}(c)=1$, there exists $r_0 \in pRp$ such that \[c + r_0(p-\tilde{s}c) \in \operatorname{U}(pRp)\]
	Label $r = r_0\theta$  and consider the element $M + r(1-sM).$
	 Since $c = pc,$  $\,\, \eta = p\eta$ and $\beta = (1-p)\beta$ we can write 
	\begin{align*}
		M + r(1-sM) &= M +r - rs\bigl(\underbrace{c\theta + \eta + \beta}_{\scriptscriptstyle M} \bigr) \\
		&= M+r - \bigl(r_{0} \underbrace{\theta s(p}_{\tilde{s}}c)\theta + r_0\underbrace{\theta s(p}_{\tilde{s}} \eta) + r_0\theta s((1-p)\beta)\bigr) \\
		&= M + r -  \bigl(r_{0}\tilde{s}c \theta + r_0\tilde{s} \eta + r_0\theta s((1-p)\beta)\bigr) \\
		&= (c\theta + \eta + \beta )+ (r_0\theta) -  \bigl(r_{0}\tilde{s}c \theta + r_0\tilde{s} \eta + r_0\theta s((1-p)\beta)\bigr) \\
		&= (c+r_0 -r_0\tilde{s}c)\theta + \left[\eta - r_0\tilde{s} \eta - r_0\theta s((1-p)\beta)\right] + \beta
	\end{align*}
Label $\eta_1 = \eta - r_0\tilde{s} \eta - r_0\theta s((1-p)\beta).$  Then $\eta_1 \in pR(1-e)$.\\
	Also label $c_1 = c+r_0 -r_0\tilde{s}c.$  Then note that since  $r_0 =r_0 p \in pRp$ we have  $ c_1 = c+r_0 -r_0\tilde{s}c = c +r_0(p -\tilde{s}c)  \in U(pRp)$  \\
	Now $c_1$ being invertible in $pRp$ and $(\theta, \theta')$ being an invertible block in $pRe$ $\implies$ $(c_1\theta, \theta'c_1^{-1}) $ is an invertible block in $pRe$ since 
\[	(\theta'c_1^{-1})(c_1 \theta )= \theta' \theta = e \,\, \text{ and } (c_1 \theta )	(\theta'c_1^{-1}) =c_1(\theta \theta')c_1^{-1} =p \]
		Then we have $M + r(1-sM) = c_1 \theta + \eta_1 + \beta$ where $c_1\theta $ is an invertible block in $pRe$ and $\beta $ is an invertible block in $(1-p)R(1-e)$ and $\eta_1 \in  pR(1-e)$. Then by Lemma \ref{triangular} $(ii)$ 
		\[M + r(1-sM) = c_1 \theta + \eta_1 + \beta \in U(R)\]
		Since $s$ is arbitrary,  by Lemma \ref{stable} we have $\operatorname{sr}_{\scriptscriptstyle R}(M) = \operatorname{sr}_{\scriptscriptstyle R}{\left(c\theta +\eta +\beta\right) }=1$ for any $\eta \in pR(1-e)$.\\
		
\noindent 		\textbf{Conversely}, let  $\operatorname{sr}_{\scriptscriptstyle R}{\left(c\theta +\eta +\beta\right) }=1$ for every $\eta \in pR(1-e)$.\\ 
Since $(\theta' + \beta')$ is invertible using Block Lemma, we can write 

\[	(c\theta +\eta +  \beta )(\theta' + \beta') = c + \eta \beta' + (1-p) = (1+ \eta \beta')(c+ 1-p) \] 
and since,  multiplication by units does not change the stable range of an element we have  $\operatorname{sr}(c + (1-p)) = 1$. Then by Lemma~\ref{stable}, for each $s\in R$ there exists $r \in R$ such that  \[ (c+ 1-p) +r (1 - s(c+1-p)) \in U(R) \]
\noindent  Fix arbitrary $s \in pRp$. We prove that there exists $r \in pRp$ such that $c + r(p-sc) \in U(pRp)$. \\
Label $s'=s+(1-p).$ Then there exists $r_1 \in R$ such that   \[V = (c+ 1-p) +r_1 (1 - s'(c+1-p)) \in U(R)\]
Observe that  $s'\bigl(c+(1-p)\bigr)=sc+(1-p) \implies 1 -s'\bigl(c+(1-p)\bigr)=p-sc \in pRp$.
\begin{align*}
\text{Then, } V \; &= \; c \, + \, (1-p) \, + \, r_1(p-sc) \\ &=  \; c \, + \, (1-p) \, + \, (pr_1 + r_1 \,- \, pr_1)(p-sc) \\
&=\;\underbrace{c\,+\, (pr_1p)(p-sc)}_{\in\,pRp} \;+\;\underbrace{(1-p)r_1(p-sc)}_{\in\,(1-p)Rp}\;+\;(1-p)\;\in\;U(R).
\end{align*}

\noindent Since $(1-p)$ is invertible in $(1-p)R(1-p)$, therefore by the converse of Prop \ref{converse}, $c+ (pr_1p)(p-sc)\in U(pRp).$  \\

\noindent Hence by Lemma~\ref{stable} we have  $\operatorname{sr}_{pRp}(c)=1$.

\bigskip
\noindent $(1)(ii)$ The proof for the second part comes directly from  $(i)$. \\
Label $p_1 = e$ and $e_1  = p$,  $\theta_1 = \theta'$ and ${\theta_1}' = \theta,$  $\beta_1 = \beta'$ and ${\beta_1}' = \beta$.\\
Then $(\theta_1, {\theta_1}' )$ is an invertible block in $p_1Re_1$ and $(\beta_1,{\beta_1}' )$ is an invertible block in $(1-p_1)R(1-e_1)$. Choose $c_1 = a \in eRe = p_1Rp_1$.\\
Then for any $\eta_1 = \tilde{\eta}  \in  p_1R(1-e_1)$ by $(1)(i)$ we get 
\[	\operatorname{sr}_{\scriptscriptstyle p_1Rp_1}(c_1)=1\iff \operatorname{sr}_{\scriptscriptstyle R}\bigl(c_1\theta_1 + \eta_1 +\beta_1\bigr)=1 \]
In other words, for $\tilde{\eta}  \in  eR(1-p) $ we get 
\[	\operatorname{sr}_{\scriptscriptstyle eRe}(a)=1\iff \operatorname{sr}_{\scriptscriptstyle R}\bigl(a\theta' +\tilde{\eta} +\beta'\bigr)=1 \]

\medskip 
\noindent $(2)(i) \text{ and } (ii):$  The proof for $(2) (i) $ can be obtained by  taking any element $c \in pRp$ and considering the element $a = \theta' c \theta \in eRe.$ and following $(1) (ii)$.  \\
\medskip\\
\noindent Additionally in each of these subcases, it is straightforward to see that $\operatorname{sr}_{\scriptscriptstyle R}(a) = 1 \,( \text{or, } \operatorname{sr}_{\scriptscriptstyle R}(c)= 1)$ \\
\end{proof}

\noindent As a  consequence we recover the Suspension Theorem in \cite{khu2}. \\

 \begin{theorem}[\cite{khu2}, Theorem 6.2]
\noindent  Let $e, f \in \operatorname{idem}(R)$ be such that $e+f=1$, and let $a \in eRe$,
$\gamma \in fRe$. Then $\operatorname{sr}_{eRe}(a) = 1 \iff \operatorname{sr}_{R}(a + \gamma + f ) = 1$.
\end{theorem}
\begin{proof}
	We can write $a + \gamma + f = (a+f)(1+\gamma)$ where $(1+\gamma) \in \operatorname{U}(R)$. Then using  idempotents $e=p,$  $\theta\ = e, \eta = 0$   and the identity block $\beta = f,$ Theorem \ref{suspension} $(2) (ii)$ above  gives us the equivalence. 
\end{proof}

\begin{theorem}[\cite{khu2},Theorem 6.6]
	Let $e, f \in \operatorname{idem}(R)$ with $e + f = 1$.
	\begin{enumerate}
		\item If $a \in eRe$, $\gamma \in fRe$ and $b \in fRf$ are such that $\operatorname{sr}_{\scriptscriptstyle eRe}(a) = 1$ and $\operatorname{sr}_{\scriptscriptstyle fRf}(b) = 1$, then
		\[\operatorname{sr}_{\scriptscriptstyle R}(a + \gamma + b) = 1\]
		\item 
		If $a \in eRe$, $\gamma \in fRe$ and $u \in U(fRf)$, then \[\operatorname{sr}_{\scriptscriptstyle eRe}(a) = 1 \iff \operatorname{sr}_{\scriptscriptstyle R}(a + \gamma + u) = 1.\]
	\end{enumerate}
\end{theorem}
\begin{proof}
$(1)$  Since $\operatorname{sr}_{eRe}(a)=1$, Theorem 3.12 gives $\operatorname{sr}_R(a+\gamma+f)=1$. Also using Theorem \ref{suspension} with $\eta = 0$ gives $\operatorname{sr}_{R}(e+b) = 1.$ Now using  the expression $(a + \gamma + f ) (e + b) = ae + \gamma e + f b= a + \gamma + b,$ and that products of stable range one elements are also stable range one \cite[Lemma ~17]{chen}, the proof follows. \\
$(2)$ Straightforward.
\end{proof}
\medskip 
\noindent Using the  Suspension Theorem, we provide an analogue to Ara and Goodearl's result \cite[Lemma ~2]{goodearl1} on stable range of corners through subequivalence.
\medskip 
\begin{cor}
	Let $e,f \in \operatorname{Idem}(R)$ such that $ e \lesssim f$. \\ If $\operatorname{sr}(fRf)=1$,
	then $\operatorname{sr}(eRe)=1$.
\end{cor}
\begin{proof}
	 There exists elements $r \in eRf$ and $s \in fRe$ such that \[(r)(s )=e\]
	Label $S = fRf$ with identity $f$ and take $g = (s)(r) \in S.$ Then clearly, $g^2 = g$ and $gf  = fg = g$. By the definition of $g$, we also have $r = r(g) \in Rg$ and $s = (g)s \in gR$. \\
	Now consider the corner ring $gSg$ in the ambient ring $S$. Define $h = f-g$. Then $h \in \operatorname{Idem}(S)$ and $hg = gh = 0$. \\
	
	\noindent We prove $\operatorname{sr}(eRe)=1$. Suppose $xa + b = e$ for arbitrary $a,b,x \in eRe.$ \\
	
	Consider the following elements in $gSg$.
	\[a' \, = \, (s)a(r) \;\;\qquad \; b' = (s)b(r) \;\;\;\; \qquad  x' = (s)x(r)\]
	\[\text{Then, } x'a'+b'= (s)(xa+b)(r) = ser = sr= g \tag{$\star$} \]
	Consider the element $a' + h \in S$. 	Then $\operatorname{sr}(S) = 1 \implies \operatorname{sr}_S(a'+h) = 1.$ 
	Since $g$ is the identity  in $gSg$ and $f-g$ is the identity in  $(f-g)S(f-g)$, by    Theorem \ref{suspension}, we have 
	\[ \operatorname{sr}_{\scriptscriptstyle gSg} (a')= 1 \iff \operatorname{sr}_S(a' + h) = 1\]
	Thus $\operatorname{sr}_{\scriptscriptstyle gSg} (a')= 1 \implies \emph{e.q.}(\star)$ is reducible. Thus there exists $y \in gSg$ such that $a' + yb' \in U(gSg)$. Label $u_g = a' + yb' \in U(gSg)$. \\
	 Then, since $y,b' \in gSg$ and $g = sr$  we have 
	 \[r(a' + yb' )s = ra's + (ry)(g)(b's) = (ra's)+ (rys)(rb's)= a+ (rys)b\]
	where $r(a' + yb' )s  = ru_gs \in U(eRe)$ since  
	\[\underbrace{r(u_g)s} \underbrace{r(u_g^{-1})s}=r(u_g)g(u_g^{-1})s = rgs=e =r(u_g^{-1})g(u_g)s = \underbrace{r(u_g^{-1})s}\underbrace{ r(u_g)s}\] 
	Thus $a+ (rys)b \in U(eRe) \implies \operatorname{sr}(eRe)=1$.

\end{proof}
\begin{remark}
A rather neat  consequence is rediscovered through  the above corollary. Suppose $e$ is a full idempotent in a ring $R$ such that $eRe$ has stable range one.  Then  $e$ being full, by Ara-Goodearl \cite[p.~3]{goodearl1}, for some $n\in\mathbb{N}$ we have
\[
E_{11}\;\lesssim\;\hat{e}=\operatorname{diag}(e,\ldots,e)\in M_n(R).
\]
Now $\hat{e}\,M_n(R)\,\hat{e}=M_n(eRe)$ has stable range one, since $eRe$ does, due to Vaserstein.
Then, applying the above corollary, we have   $\operatorname{sr}(R) = 1$. 
\end{remark}

\medskip 

\section{Block Invertibility and Special-Clean elements}

\medskip 
Abrams and Rangaswamy \cite{ranga} coined the term \emph{special clean} for
rings whose elements decompose in the following way :\\
\emph{ An element $a\in R$ that can be written in the form $a=e+u$ with $e^2=e$, $u\in U(R)$ and $aR \, \cap \, eR=\{0\}$}.\\
They used it to compare regularity conditions for Leavitt path algebras. Camillo and Khurana \cite{cam} proved that on the entire ambient ring, unit regularity is equivalent to special clean. However, locally a unit regular element is not necessary special-clean (or even clean), and there are numerous such examples, one of which was provided   in \cite{khulam} where Khurana and Lam  showed the element $\bigl(\begin{smallmatrix} 1+xy & x^2 \\ 0 &0 \end{smallmatrix} \bigr)$ is unit-regular but not clean over $K= k[x,y]$ for any commutative domain $k$. Many authors  investigated additional conditions under which unit-regular elements are special clean \cite{ster, khu4,patty,zhou1}. Nielsen and Ster gave an example in \cite{ster1} where they showed  that even when $a, 1-a$ is unit regular, even then $a$  need not be clean. They also used regularity of the element  $a$ say $a= ara$ and unit regularity of the specific corner element $eae$ where $e = ra,$  to obtain a stronger condition than special clean (see below). \\

\begin{theorem}{\cite[Theorem~3.14]{ster1}}
Let $a$ be an element of a ring $R$. The following are equivalent:
\begin{enumerate}
\item The exists $r\in R$ with $ara= a$, and writing $e:= ra\in \operatorname{idem}(R)$, the element $eae$ is unit-regular in $eRe$.
\item There exist $e\in \operatorname{idem}(R)$ and $u\in U(R)$ such that $a= e+u$, $aR\cap eR= (0)$, and $a^2R\cap aeR= (0)$.
\end{enumerate}
\end{theorem}

\medskip 
Mary and Patricio identified where stable range can be applied to evoke special-cleaness of $a$.  In particular, they obtained that when $av^{-1}a = a$, and the Pierce decomposition of $v$  with respect to $f = av^{-1}$ is given by : $v  \equiv \bigl(\begin{smallmatrix} v_1& v_2 \\ v_3 & v_4\end{smallmatrix} \bigr)$  and if there exists a $\gamma \in \bar{f}Rf$ (resp. $\sigma \in fR\bar{f}$ ) such that  $v_4 + \gamma v_2 \in U(\bar{f}R\bar{f})$  (resp. $v_4 + v_3\sigma \in U(\bar{f}R\bar{f}$) then $a$ is special clean. Subsequently, they proved the following: 
\medskip 

\begin{cor}{\cite[Cor~3.4]{patty}}
	Let $av^{-1}a= a$ and $f = av^{-1}$. If $\operatorname{sr}(v_4)=1$ in $ \bar{f}R\bar{f}$ (in particular if $v_4$ is unit-regular in $\bar{f}R\bar{f}$), then $a$ is special clean.
\end{cor}

We give an element-wise characterization under which an unit-regular element is special-clean. Recall Bass's stable range condition as follows: A ring $R$ has stable range 1 if \underline{every} left unimodular pair \((a, b)\) is reducible, that is whenever $ra + sb = 1 \implies  (a + xb)$ is invertible for some \(x \in R\) for arbitrary $a,b \in R$.  On a specific element $a$, $\operatorname{sr}(a) = 1$ is defined as every left (right) unimodular pair of the form $(a,y)$ is left (right) reducible for arbitrary $y \in R$. We observed that the missing link is simply reducibility  for a  specific unimodular pair $(u_4, -)$ where $a = aua$ and the Pierce decomposition of $u$ is \[ u = \underbrace{(au)u(au)}_{u_1} + \underbrace{(au) u(1-au)}_{u_2} + \underbrace{(1-au) u(au)}_{u_3} + \underbrace{(1-au)u (1-au)}_{u_4}\]

We show below that a necessary and sufficient condition for $a$ to be special clean requires that  the specific  pair $(u_4, \psi)$ where $u_4 = (1-au)u(1-au), \,\, \psi = (1-au)(u^{-1}au)u(1-au)$ which is left unimodular in the corner ring $(1-au)R(1-au)$  is left reducible: \emph{i.e.} $\alpha+z\psi\in U(S)$ for some $z\in S$. \\

\noindent We first consider the following lemma for any element $a = aua:$

\begin{lemma}{\label{piercedecom}}
	Let $ a\in R $ such that $a = aua$ with $u  \equiv \bigl(\begin{smallmatrix} u_1& u_2 \\ u_3 & u_4\end{smallmatrix} \bigr)$ as above. Then for any fixed $ \sigma \in (1-au)R(au)$  the following are equivalent:
	\begin{enumerate}[label=(\roman*)]
		\item  $u(1-au) + (au )+ \sigma \in U(R)$
	\item  $u_4 - \sigma u_2 \in U\bigl((1-au)R(1-au)\bigr)$
\end{enumerate}
\end{lemma}
\begin{proof}
	Observe that $u_4 - \sigma u_2 = (1-au)u(1-au) - \sigma u(1-au) $.\\
	
	\noindent $(i) \iff (ii):$ Let  $ \sigma \in (1-au)R(au)$ be such that  \[u(1-au) + (au )+ \sigma \in U(R)\] 
	
	\noindent Write $p = au$. Label $\alpha = (1-p)u(1-p)$ and $\gamma = pu(1-p) \in pR(1-p).$ \\
	Label $V = u(1-au) + (au) +  \sigma = (\alpha + \gamma )+ p + \sigma.$\\
	Now $\sigma \in (1-p)Rp \implies (1-p)\sigma p = \sigma.$ Similarly, for $\gamma$, we have $p\gamma(1-p) = \gamma.$ Then clearly,  $\sigma^2=0=\gamma^2 \implies (1+\sigma),\,(1+\gamma)\in U(R)$. \\
	
	\noindent Then, expanding the product $ (1+ \sigma) (\alpha-\sigma\gamma+p)(1+\gamma)$ we get 
	\begin{align*}
		(1+\sigma)\bigl[(\alpha-\sigma\gamma+p)(1+\gamma)\bigr]
		&=(1+\sigma)\Bigl[\alpha-\sigma\gamma + p
		+\underbrace{(\alpha-\sigma\gamma)\gamma}_{0}
		+\underbrace{p\gamma}_{\gamma}\Bigr]\\
		&=\alpha-\sigma\gamma+p+\gamma
		+\underbrace{\sigma(\alpha-\sigma\gamma)}_{0}
		+\underbrace{\sigma p}_{\sigma}+\sigma\gamma \\
		&=\alpha+\sigma+\gamma+p
	\end{align*}
	
	\noindent  Then $1+\sigma$ being invertible shows that 
	\[
	V\in U(R)\iff\underbrace{\alpha-\sigma\gamma}_{(1-p)R(1-p)} \,\,\, + \,\,\, p \,\,\, \in U(R)
	\]
	
	\noindent Since $ p \in U(pRp)$  and $\alpha-\sigma\gamma+p \in U(R)$, then by Proposition \ref{converse},  $ \alpha-\sigma\gamma$ must have a block inverse in $(1-p)R(1-p)$. in other words, 
	\[u_4 - \sigma u_2= \alpha - \sigma \gamma  \in U\bigl((1-p)R(1-p)\bigr)\]

	\noindent The above line is actually  reversible using Lemma \ref{triangular} and Proposition \ref{converse} that is, 
	\[ (u_4 - \sigma u_2) \,\,\, +\,\,\, p \,\, \in U(R)
	\iff (u_4 - \sigma u_2)\in U\bigl((1-p)R(1-p)\bigr) \]
	Thus the converse also holds.  This proves (i)$\iff(ii)$, with the same
	pair $u$ and $\sigma$ on both sides. \\
\end{proof}

For our construction of clean elements, we follow the strategy given  in ({Corollary 27,\cite{siddique1}}) and (\cite{mary2} Corollary 2.3.), namely  for an element $a \in R$, $a$ is clean if and only if there exists outer inverses $x, y $ of $a$ and $(1-a)$ resp.  such that $ax+ (1-a)y= 1$ and $(x-y)$ is invertible.

	\begin{theorem}{\label{special-clean}}
	Let $ a\in R.$ Then the following are equivalent :
	\begin{enumerate}[label =(\roman*)]
		\item $a = f+ v$ for  $f^2 = f$ and unit $v \in U(R)$ such that $aR \cap fR = \{0\}$\\
			\item There exists  unit $u \in U(R)$ and $ \sigma \in (1-au)R(au)$ such that $a = aua$ and
		\[(1-au)u(1-au) - \sigma u(1-au) \in U\bigl((1-au)R(1-au)\bigr)\]
		\item there exists $u\in U(R)$ with $a=aua$ such that the pair
	$(\alpha,\psi)$, where
	\[
	\alpha=(1-au)u(1-au),\qquad \psi=\mu\gamma=(1-au)(u^{-1}au)\,u(1-au),
	\]
	which is left unimodular in $S=(1-au)R(1-au)$, is left reducible:
	$\alpha+z\psi\in U(S)$ for some $z\in S$.
	\end{enumerate}

	\end{theorem}
	\begin{proof}
		$(ii) \implies (i):$ 
Let there  be a unit $u \in U(R)$ and $ \sigma \in (1-au)R(au)$ such that $a = aua$ and \[u(1-au) + (au )+ \sigma \in U(R)\]
We will find  outer inverses $x,y$ of $a$ and $(1-a)$ resp. with $ax +(1-a)y = 1,$ such that $(x-y)$ is invertible and $aR\, \cap \,\,(x-y)(1-ax)R = \{0\}$. \\

\noindent Define $p = au,\,\,\,\,\, q = ua, \,\,\,\, \, e = p +a(1-p), \,\,\,\ \text{and } \, y = 1-p $
 
 \medskip 
 \noindent Then clearly  $p^2=p, \, $ $\, q^2=q, \, $ $\, pa=a; \, $
  \setcounter{claim}{0} 
\begin{claim}
	$e^2 = e; \,\,ea = a;$  $\,\, pe = e; \,\, $	\( \left\{
	\begin{aligned}
		(1-p)(1-e) &=1-p\\
		(1-e)(1-p) &=1-e
	\end{aligned}
	\right.
	\) 
\end{claim}
\begin{proof}
Using $pa=a$ and $(1-p)a=0$ we get
\begin{align*}
	e^2 &= \left(p +a(1-p)\right) \left(p +a(1-p)\right) \\
	&= p^2+pa(1-p)+a(1-p)p+a(1-p)a(1-p) \\
	&=  p + a(1-p) = p + (e-p) = e
\end{align*}

\noindent Using $pa = a$ we also get 
\[ea=\left(p +a(1-p)\right)a = pa+a(1-p)a=a\]
\noindent Also, using  $pa = a$ we get $pe=p+a(1-p)=e$ and therefore,
\[(1-p)e = 0 \Longleftrightarrow (1-p)(1-e) = 1-p\]
Finally we have $ep = p \Longleftrightarrow (1-e)p = 0 \Longleftrightarrow (1-e)(1-p) = 1-e$\\
\end{proof}
\noindent Now label $y = (1-p)$. Then clearly, $ya=(1-p)a=0$. \\
Therefore, $y(1-a)y = y^2 = y \implies $ y is an outer inverse of $(1-a)$. Finally, using $e = p + a(1-p)$ gives $(1-e) = (1-p) - a (1-p) = (1-a)(1-p) = (1-a)y$. Thus \[e +(1-a)y = 1\]

\medskip
 \noindent We still need to obtain the outer inverse $x$ of $a$ such that $ax = e$ which would imply that $ax + (1-a)y = 1$. \\
 
 \noindent Define
 \[
 \alpha = (1-p)u(1-p), \quad  \gamma = pu(1-p), \quad \lambda = (1-p)u^{-1}(1-p), \quad \mu = (1-p)u^{-1}p \]
 Now by assumption,  there exists $\sigma \in (1-p)Rp$ such that 
 \[(1-au)u(1-au) - \sigma u(1-au) \in U\bigl((1-au)R(1-au)\bigr)\]
 Then for the same $\sigma$ and $u$ by Lemma \ref{piercedecom} we have an equivalent statement
 \[u(1-au) + (au )+ \sigma = \alpha + \gamma + p + \sigma \in U(R)\]
 
 \noindent Label $V = \alpha + \gamma + p + \sigma$.  Post-multiplying $V$ with $p$ and $(1-p)$ respectively gives 
 \[	 (V)(p) = \alpha p + \gamma p + \sigma p + p = \sigma + p \]

\noindent and	$(V)(1-p) = \alpha + \gamma = u(1-au)= (1-ua)u \hspace{19mm} \cdots \hspace{16mm} (B^{\star})$ \\

	\noindent 			Now label $Vp  =  \sigma + p = g_1$ and $\pi = VpV^{-1}$. Then $g_1, \pi  \in \operatorname{Idem}(R)$ as seen below.
	\begin{align*}
	{g_1}^2 &= (\sigma + p)^2 \hspace{10mm } \cdots \hspace{5mm} \text{ where } \sigma \in (1-p)Rp\\
	&= \underbrace{{\sigma}^2}_{0} + \underbrace{\sigma p}_{\sigma} + \underbrace{p\sigma}_{0} + \,\, p = \sigma + p = g_1
\end{align*}
\item And, $\pi^2 = \left(VpV^{-1}\right)\left(VpV^{-1}\right)= VpV^{-1}= \pi$. \\
We also have $\pi g_1= (VpV^{-1})(Vp)= Vp = g_1 $ and   \[g_1 \pi = \left(Vp\right) \left(VpV^{-1}\right)= \underbrace{\left(Vp\right)}_{g_1}\underbrace{\left(Vp\right)}_{g_1}V^{-1}= \left(Vp\right)V^{-1}= \pi \]
From $(B^\star),$ $(1-\pi) = V(1-p)V^{-1} = (1-ua)uV^{-1} \implies (1-ua)(1-\pi)= 1-\pi \implies a\pi = a.$ 
Furthermore, from $(B^{\star})$ we also get 
\begin{align*}
(1-\pi)(1-ua) & =V(1-p)V^{-1} V(1-p)u^{-1} = V(1-p)u^{-1} = (1-ua) \\
\implies   \pi (ua) &= \pi \hspace{40mm} \cdots \hspace{30mm} (C^{\star})
\end{align*}

\medskip 
Now choose $ x = \pi ue \in Re$. Then $ax = pe=e$. Also $x(ax) = xe = x$ and $axa = ea = a$ from Claim 1. Hence $x \in Ref(a)$. \\
Set $d = x-y$. Then since $(1-p)(1-e)= (1-p)$ we can write \[d = \underbrace{px}_{\in pRe} \,\, \, +  \,\, \, \underbrace{(1-p)x}_{(1-p)Re} \,\,\ + \,\,\,  \underbrace{\left(-(1-p) \right)}_{(1-p)R(1-e)}\]

\noindent Set $\beta = -y = -(1-p) $ and $ \beta' = -(1-e) \in (1-e)R(1-p).$ Then using Claim 1, we get $\beta  \in (1-p)R(1-e) $ and $\beta' \in (1-e)R(1-p)$ such that \[\beta \beta' = (1-p)(1-e) = (1-p) \,\,\, \text{ and } \,\,\, \beta' \beta = (1-e)(1-p) = (1-e)\] 

\noindent Also, with $ea = a, \,$ and $g_1p = g_1,$ for  elements $ag_1 \in eRp$,  and $px \in pRe $ we obtain
\[
(px)(ag_1) = p(xa)g_1 = p(\pi u e a)g_1 = p(\underbrace{\pi ua}_{C^{\star}})g_1 =  p(\pi)g_1 = pg_1 = p
\]
\[\text{and, } \, \,\,(ag_1)(px) = a(g_1p)x = ag_1x = ag_1 (\pi u e )= a (g_1 \pi ) ue =a\pi ue = ax=e\] 
\medskip 

\noindent Thus, $px$ and $\beta$ are invertible blocks in $pRe$ and $(1-p)R(1-e)$ with inverses $ag_1$ and $\beta'$.
\noindent Therefore, by our Block Lemma \ref{triangular},  $d$ is invertible in $R$. \\
In fact, Lemma \ref{triangular} gives us its inverse to be the sum of the block inverse of $px$ in $eRp$, block inverse of $\beta,$ and  $(\beta) ((1-p)x)(ag_1) $ \[d^{-1} = ag_1 + \beta'  - \beta' (1-p)x(ag_1) \]

\noindent Thus, $x-y$ being invertible proves that $a$ is clean. For the sake of completition we show it below. \\
$	ad - d(1-e) = a(x-y) - (x-y)(1-e) = ax - ay - \underbrace{x(1-e)}_{0} + \underbrace{y(1-e)}_{y} = ax - ay + y $ \\
$\implies ad - d(1-e)=e + (1-e) = 1 \implies a = d^{-1} + d(1-e)d^{-1}  $

\noindent We  observe that $a$ is actually special-clean. \\
Set $f= d(1-e)d^{-1}$. Then $a - d^{-1} = f$. \\
Suppose $ \exists r \in aR \cap fR $. Then we can write $r = a r_1 = fs_1$ for $r_1,s_1 \in R$. Now 
\[
ada=axa-a(ya)=a-0=a ,
\]
Also  $f=a-d^{-1} \implies adf = ad(a-d^{-1})=ada-ad\,d^{-1}=a-a=0 $

\noindent Then $adr = adar_1 = ar_1 = r$ while $adr = adfs_1 = 0$. Hence $r = 0.$ \\
Thus $aR \cap fR = \{0\}$ and $a$ is therefore, special-clean.\\

\bigskip 
\noindent $(i) \implies (ii):$ Suppose there exists  $v \in U(R)$ and $f^2 = f$ such that $a = f+v$ where $aR \cap fR = \{0\}$.  Take $ u = v^{-1}$ and $ p = au$. \\
Then $a u = fu + 1 \implies (au-1) = fu \in fR \implies (au-1)a \in fR$ \\
Together with $(au-1)a = a(ua-1) \in aR	$, we must have $(au-1)a = 0$ since $aR \cap fR = \{0\}$. Thus $a = aua$. \\
Then clearly $p^2 = p$ and  $ p(f )= au(a-v) = aua - a = 0$. Consequently, $(1-p)f = f$ and $p(1-f)= p$. \\
Also from above $ (au-1) = fu  \implies (1-f)(p-1) = 0 \implies (1-f)p = 1-f,$  and $ f(1-p) = (1-p).$\\

\noindent We are now ready to construct our $\sigma \in (1-p)Rp$. \\
Set $ \sigma = (1-f) - p.$  Then using $(1-f)p = 1-f$ gives  \[\sigma = -fp = -(1-p)fp \in (1-p)Rp\]
Label $V = u(1-au) + (au) +  \sigma  = u(1-p) + p + \sigma = u(1-p) +p+ (1-f-p) = u(1-p) + (1-f).$ We prove that $ V \in U(R)$ using  Block Lemma.\\

\noindent Using   $fu = -(1-p), f(1-p) = 1-p $ and $(1-f)p = 1-f$,  we  rewrite $V$ as 
\begin{align*}
V= u(1-p) + (1-f) &= (1-f)u(1-p) + fu(1-p) + (1-f)p \\
&= (1-f)u(1-p) -(1-p) + (1-f)p \\
&= (1-f)u(1-p) -f(1-p) + (1-f)p \\
\end{align*}

\noindent Thus we can write  $ V = \underbrace{(1-f)p}_{(1-f)Rp}  \,\, + \,\, \underbrace{(1-f)u(1-p)}_{(1-f)R(1-p)} \,\ - \,\  \underbrace{f(1-p)}_{fR(1-p)}$. To use the Block lemma, we show the extreme components are  inverses in $(1-f)Rp$ and $fR(1-p).$\\

\noindent Using $p(1-f) = p $  and $(1-f)p = 1-f $ we get \[\underbrace{(1-f)p}_{\in (1-f)Rp} \,\cdot \, \underbrace{p(1-f)}_{\in pR(1-f)}  = (1-f)  \,\,\,\,\, \text{ and }  \,\,\,\,\,  \underbrace{p(1-f)}_{\in pR(1-f)} \,\, \underbrace{(1-f)p}_{\in (1-f)Rp}  = p\].

\noindent and using $f(1-p) = (1-p)$ and $(1-p)f = f$ we get \[ \underbrace{f(1-p)}_{\in fR(1-p)} \,\,\,\, \underbrace{(1-p)f}_{\in (1-p)Rf}  = f  \,\,\,\,\, \text{ and }  \,\,\,\,\,  \underbrace{(1-p)f}_{\in fR(1-p)} \,\,\,\, \underbrace{f(1-p)}_{\in (1-p)Rf} = 1-p\]

\noindent Thus $(1-f)p$ is invertible in $(1-f)Rp$ and $f(1-p)$ is invertible in $fR(1-p)$. Then by Block Lemma \ref{triangular} $(ii),$ $V$ is  invertible in $R$\\ 
Hence by Lemma \ref{piercedecom}, there exists a $\sigma \in (1-au)R(au)$ such that 
\[(1-au)u(1-au) - \sigma u(1-au) \in U((1-au)R(1-au))\]

\medskip 
\medskip

\noindent $(iii) \implies (ii):$ Label  $S = (1-au)R(1-au)$. Let  $u \in U(R)$ such that $a = aua $ and $(\alpha, \psi)$ is unimodular in $S$ where  $\alpha=(1-au)u(1-au)$ and $\psi = (1-au)(u^{-1}(au)u)(1-au)$.\\
\noindent Since $(\alpha, \psi)$ is left reducible in $S$, there exists  $z\in S$ such that\[ \alpha+z\psi = \alpha + z(\mu \gamma)\in U(S)\]
Label $\sigma = - z((1-au)u^{-1}(au))  \in (1-au)R(au).$  Then $\sigma (u)(1-au) = -z \psi$ and 
\[
(1-au)u(1-au)-\sigma u(1-au) = \alpha+z\psi \in U(S).
\]

\medskip 

\medskip
\noindent $(i) \implies (iii):$ 	If $a$ is special clean, then there exists $f^2 = f$ and $v \in U(R)$ such that  \[a=f+v  \,\,\,\ \text{ and } \,\,\, aR\cap fR=\{0\} \] 
Take $u=v^{-1}$ and $p=au$.  Consider the pair $(\alpha,\psi)$, where  $\alpha=(1-au)u(1-au)$ and $ \psi =(1-au)(u^{-1}au)\,u(1-au)$. It can easily be seen that this pair is  left unimodular in $S=(1-au)R(1-au).$

\noindent Then from   the proof in $(i)\iff (ii),$  for $\sigma = 1-f-p = -fp\in (1-p)Rp $ we have \[ \alpha-\sigma u(1-au)\in U\bigl((1-p)R(1-p)\bigr)  \hspace{3mm} \cdots \hspace{3mm} (\star\star)\] 

\noindent Now  $(1-p)u^{-1}p=(1-p)vp=(1-p)(a-f)p
=\underbrace{(1-p)a}_{0}p-(1-p)fp$ where  
\begin{align*}
	(1-p)fp &= (1-au)(a-u^{-1})(p) \\
	&= (1-au)(-u^{-1})(p)\\
	&= (a - u^{-1} )(p) = fp
\end{align*}
Therefore, $(1-p)u^{-1}p =-fp=1-f-p=\sigma$. Then $\psi = (1-au)(u^{-1}au)\,u(1-au) = \sigma u(1-au)$. Then from $( \star\star)$ above, \[\alpha  + (-1_{\scriptscriptstyle S}) \psi \in U(S)\]

\noindent Hence  the pair $(\alpha,\psi)$ is left reduced in $S$.

\end{proof}

\begin{remark}
Applying this theorem on the entire ring enhances a known characterization of unit-regular rings. It is quite well known that  unit-regularity is encoded by a canonical internal decomposition $R=aR\oplus (a-u)R$ and can be found in several papers such as \cite{calgu}, \cite{tklee1}, \cite{zhou3}.  We can further enhance and  construct what the element $a-u$ looks like : \\
\emph{A ring $R$ is unit regular if and only if for every element $a$ there exists unit $u \in U(R)$ and $f \in  \operatorname{Idem}(R)$ such that $a = aua $ and $R = aR \oplus fR$ where   $ (1-au)  - f  \in (1-au)R(au)$.}
\end{remark}
\medskip 

As a consequence of the above theorem, we recover and strengthen Mary and Patricio's \cite{patty} result to show that it is an equivalent characterization than a consequence. The proof   shows that the $y v_2 \,\, + \,\, v_4$ defined below is the inverse of the element $(1-au)u(1-au) - \sigma u(1-au)$ in  $(ii)$ of the above theorem.

\begin{lemma}{\cite[Lemma~3.1]{patty}}
	Let $a = aua$ and $f = au$. Consider the Peirce decomposition  of $u^{-1}$ with respect to $(au)$ as 
	 \[ u^{-1} = \underbrace{(au)u^{-1}(au)}_{v_1} \,\,+ \,\,  \underbrace{(au) u^{-1}(1-au)}_{v_2} \,\, + \,\, \underbrace{(1-au) u^{-1}(au)}_{v_3} \,\, + \,\,  \underbrace{(1-au)u^{-1} (1-au)}_{v_4}\]
	
	If $y v_2 \,\, + \,\, v_4 \,($resp. $v_3 x \,\, + \,\, v_4)$ is a unit in $\bar{f}R\bar{f}$ for some $y \in \bar{f}R f \,\, ($resp.\ $ x \in fR \bar{f})$, then $a$ is special clean. \\
\end{lemma} 

\begin{proof}
	We prove that it is a characterization.\\
	
	Label $\bar{f} = (1-au)$.  Then   that any  $y \in \bar{f}Rf$
\[	 yv_2 \, + \, v_4 =  y (u^{-1})(1-au) \, + \, (1-au)u^{-1}(1-au)  = (y  + \bar{f})u^{-1}\bar{f}\]

\noindent Suppose  $V	=( \bar{f}+y)u^{-1} \bar{f} $ is a unit in $\bar{f}R\bar{f}$ for some $y\in \bar{f}R f$. Then \[V^{-1} V = V^{-1} ( \bar{f} \, +\, y)u^{-1} \bar{f}  = \bar{f}\]

\noindent  Label  $\sigma_0=V^{-1}(\bar{f} \,+ \, y)u^{-1} \,\, \in \bar{f} R$.
		Then $\sigma_0 \, \bar{f}=V^{-1}V=\bar{f} \implies \sigma_0 (f)=\sigma_0- \bar{f}$. \\
		Now label $\sigma =  -\sigma_0(f) = -\sigma_0 +  \bar{f} \in \bar{f}Rf$. Then  we obtain
		\[
	(1-au)u(1-au) - \sigma (u)(1-au)=	\bar{f}u\bar{f}- \sigma u \bar{f} \,=\, 	\bar{f}u	\bar{f} \, + \, (\sigma_0- 	\bar{f})u	\bar{f} \, =\, \sigma_0\,u \bar{f} \]
	Now using $\sigma_0=V^{-1}(\bar{f} \,+ \, y)u^{-1}$ and $y \bar{f}=0$  we obtain 
		\[ \sigma_0\,u \bar{f} = V^{-1}(\bar{f} \,+ \, y)u^{-1}u\bar{f}\,=\, V^{-1}\bar{f}\,=\, V^{-1}	\]
		Thus $	(1-au)u(1-au) - \sigma (u)(1-au)$ is the inverse of $yv_2 \, + \, v_4 $ for some $y \in \bar{f}R f$. 
	Then by condition $(ii) $ of the above Theorem \ref{special-clean},  $a= aua$ is special clean. \\
	
	Conversely, suppose  $a=e+v$ be special clean and put $u=v^{-1}$, so that $aua=a$
	and $1-au=-ev^{-1}$. Then $(1-au) \in \operatorname{Idem(R)}$ gives $(1-au)v=-(1-au)e$, while
	$e(1-au)=1-au$, so
	\[
	v_4=(1-au)u^{-1}(1-au)=-(1-au)e(1-au)=-(1-au),
	\]
which is a unit of $(1-au)R(1-au)$. Hence the hypothesis holds with $y=0$.
	
\end{proof}
\begin{remark}
	It is not too difficult to see that in the opposite ring $R^{\operatorname{op}}$, for the element $u^{-1}au$ where  $(u^{-1}au)^{-1} \circ u \circ  (u^{-1}au)= u^{-1}au$, the Peirce decomposition of $u^{-1}$ with respect to the idempotent $(u^{-1}au )\circ u = au$ transforms as
	\[
	v_4^{\operatorname{op}}=v_4,\,\,\,\,\,\, v_2^{\operatorname{op}}=v_3,\,\,\,\,\,\, x = y^{\operatorname{op}}\in (1-au)R^{\operatorname{op}}(au)=(au)R(1-au)
	\]
	
	so the term $y^{op}v_2^{op}+v_4^{op}\in U(S)$  becomes  $v_3x+v_4\in U(S) $. Thus $u^{-1}au$ is special clean $\iff$ $a$ is special clean. (using  symmetricity of left ideal  and right ideal intersection in a special-clean decomposition as proved in  \cite{ster} Theorem 2.3).
\end{remark}

\begin{cor}{\label{corner}}
	Let $a\in R$ with $a=aua$ for some $u\in U(R)$. Label the idempotent $ p = au$ and the corner ring $S = (1-p)R(1-p)$. If 
	\[
	\operatorname{sr}_{\scriptscriptstyle S}\bigl((1-p)u(1-p)\bigr)=1
	\]
	then $a$ is special clean.
\end{cor}

\begin{proof}

	\[
\noindent \text{Label } \,\, \,\,  	\lambda=(1-p)u^{-1}(1-p),\quad
	\mu=(1-p)u^{-1}p,\quad
 \,\,\,\  \alpha  \,= (1-p)u(1-p)
	\]
	and $\gamma = pu(1-p)$. Label $	\psi=\mu\gamma$. Then clearly $\lambda, \alpha, \psi \in (1-p)R(1-p)$ and 
	\begin{align*}
		\lambda\alpha+\psi
		&=(1-p)u^{-1}(1-p)u(1-p)+(1-p)u^{-1}pu(1-p)\\
		&=(1-p)u^{-1}\bigl((1-p)+p\bigr)u(1-p)=(1-p)=1_{\scriptscriptstyle S} 
	\end{align*}
	Since $\operatorname{sr}_{\scriptscriptstyle S}\bigl(\alpha \bigr)=1$ there exists  $z\in S$ such that 
	\[
	\alpha+z\psi=\alpha+z\mu\gamma\in U(S).
	\]
	Label $\sigma=-z\mu$. Then $\sigma\in(1-p)Rp$. Then 
\[\alpha-\sigma\gamma= (1-au)u(1-au) \,\, - \,\, \sigma u(1-au) \in U((1-au)R(1-au))\]
Then by   above Theorem \ref{special-clean} $(ii)$, $a$ is special-clean.
\end{proof}

\noindent We can  recover the following result of Camillo and Khurana on the entire ring $R$. 

\begin{theorem}{\cite{cam}}
	A ring $R$ is unit-regular if and only if every element $a$ of R can be written as $e + u$ such that $aR \cap eR = \{0\}$, where $e$ is an idempotent and $u$ a unit in $R$.
\end{theorem}
	\begin{proof}
Since $R$ is unit regular, $R$ has stable range 1 and consequently, for any idempotent $e \in R$, $eRe$ also has stable range 1. For any element $a \in R$, there exists $u \in U(R)$ such that $a = aua$. Set $e = (1-au)$. Then $\operatorname{sr}_{eRe}(eue) = 1$ and subsequently, by Corollary \ref{corner}, $a$ is special clean. \\
	The converse is straightforward.\\
	\end{proof}
	

	\end{document}